\documentclass[aap,MSNbibl,citesort,dvips]{arximspdf}
\usepackage{mathrsfs}

\doi{10.1214/12-AAP851} 
\volume{22}
\issue{6}
\pubyear{2012}
\firstpage{2505}
\lastpage{2538}

\makeatletter
\newtheorem{theorem}{Theorem}[section]
\newtheorem{lemma}[theorem]{Lemma}
\newproclaim{remark}[theorem]{Remark}

\newproclaim{definition}[theorem]{Definition}
\newtheorem{proposition}[theorem]{Proposition}
\newproclaim{Examples}[theorem]{Examples}

\def\eps{\varepsilon}
\def\div{\operatorname{div}}
\def\dif{\mathrm{d}}
\def\cB{\mathcal{B}}
\def\cG{\mathcal{G}}
\def\cL{\mathcal{L}}

\def\sA{\mathscr{A}}
\def\sF{\mathscr{F}}
\def\sG{\mathscr{G}}
\def\sK{\mathscr{K}}
\def\sN{\mathscr{N}}

\def\mB{\mathbb{B}}
\def\mD{\mathbb{D}}
\def\mE{\mathbb{E}}
\def\mI{\mathbb{I}}
\def\mN{\mathbb{N}}
\def\mR{\mathbb{R}}
\def\mW{\mathbb{W}}

\makeatother

\begin{document}
\begin{frontmatter}

\title{Stochastic functional differential equations driven by L\'evy processes
and quasi-linear partial integro-differential equations\thanksref{T1}}
\runtitle{SFDE driven by L\'evy processes and quasi-linear PIDE}
\thankstext{T1}{Supported by NSFs of China (Nos. 10971076; 10871215) and the
Program for New Century Excellent Talents in University (NCET-10-0654).}

\begin{aug}
\author[A]{\fnms{Xicheng} \snm{Zhang}\corref{}\ead[label=e1]{XichengZhang@gmail.com}}
\runauthor{X. Zhang}
\affiliation{Wuhan University}
\address[A]{School of Mathmatics and Statistics\\
Wuhan University\\
Wuhan, 430072\\
P. R. China\\
\printead{e1}} 
\end{aug}

\received{\smonth{5} \syear{2011}}
\revised{\smonth{2} \syear{2012}}

%
\begin{abstract}
In this article we study a class of stochastic functional differential
equations driven by L\'evy processes
(in particular, $\alpha$-stable processes),
and obtain the existence and uniqueness of Markov solutions in small
time intervals.
This corresponds to the local solvability to a class of quasi-linear
partial integro-differential equations.
Moreover, in the constant diffusion coefficient case, without any
assumptions on the L\'evy generator,
we also show the existence of a unique maximal weak solution for a
class of semi-linear
partial integro-differential equation systems under bounded Lipschitz
assumptions
on the coefficients. Meanwhile, in the nondegenerate case
(corresponding to $\Delta^{\alpha/2}$ with
$\alpha\in(1,2]$), based upon some gradient estimates,
the existence of global solutions is established too. In particular,
this provides a
probabilistic treatment for the nonlinear partial integro-differential
equations, such as the
multi-dimensional fractal Burgers equations and the fractal scalar
conservation law equations.
\end{abstract}

\begin{keyword}[class=AMS]
\kwd{60H30}
\kwd{35R09}
\end{keyword}
\begin{keyword}
\kwd{L\'evy processes}
\kwd{Feyman--Kac formula}
\kwd{fractal Burgers equation}
\end{keyword}

\end{frontmatter}

\section{Introduction}\label{sec1}
Consider the following multi-dimensional fractal Burgers equation in
$\mR^d$:
%
\begin{equation}\label{EL5}
\partial _t u=\nu\Delta^{\alpha/2} u-(u\cdot\nabla u),\qquad t\geq0,\  u_0=\varphi
,
\end{equation}
where $u=(u^1,\ldots,u^d)$ and $\nu>0$ is a viscosity constant, and
$\Delta^{\alpha/2}$ with $\alpha\in(0,2)$ is the usual fractional
Laplacian defined by
\[
\Delta^{\alpha/2}u(x):=\lim_{\eps\downarrow0}\int_{|z|\geq\eps
}\frac
{u(x+z)-u(x)}{|z|^{d+\alpha}}\,\dif z.
\]
This is a typical nonlinear partial integro-differential equation and
is regarded as a simplified
model for the classical Navier--Stokes equation when $\alpha=2$.
Recently, there has been great interest in studying the multi-dimen\-sional
Burgers turbulence (cf.~\cite{Be-Kh,Woy}), the fractal Burgers equation
(cf.~\cite{Bi-Fu-Wo,Ki-Na-Sh,Ch-Cz-Si}) and the fractal conservation
law equation (cf.~\cite{Dr-Im}), etc.
All these works are based on the analytic approaches, especially the
energy method,
Duhamel's formulation and the maximum principle.

The purpose of the present paper is to give a probabilistic treatment for
a large class of quasi-linear partial integro-differential equations.
Let us first introduce the main idea.
By reversing the time variable, one can write Burgers' equation (\ref
{EL5}) as the following equivalent backward form:
%
\begin{equation}\label{Es3}
\partial _t u+\nu\Delta^{\alpha/2} u-(u\cdot\nabla u)=0,\qquad t\leq0,\  u_0=\varphi.
\end{equation}
Now, consider the case of $\alpha=2$, and for a given smooth solution
$u_t(x)\in C^\infty_b(\mR^d;\mR^d)$ to the above equation,
let $X_{t,s}(x)$ solve the following stochastic differential equation
(abbreviated as SDE):
%
\begin{equation}\label{Es1}
\dif X_{t,s}(x)=-u_s(X_{t,s}(x))\,\dif s+\sqrt{2\nu}\,\dif W_s,\qquad s\in
[t,0],\  X_{t,t}(x)=x,
\end{equation}
where $(W_s)_{s\leq0}$ is a $d$-dimensional standard Brownian motion
on $\mR_-:=(-\infty,0]$.
By It\^o's formula and the Markov property of the solution, it is well
known that
%
\begin{equation}\label{Es2}
u_t(x)=\mE\varphi(X_{t,0}(x)).
\end{equation}
Conversely, assume that $(u,X)$ solves the implicit system (\ref{Es1})
and (\ref{Es2}); then
$u$ also solves the backward Burgers' equation (\ref{Es3}). This type of
implicit stochastic differential
equation has been systematically studied by Freidlin~\cite{Fre},
Chapter~5; see also~\cite{Bl,Ta}.

Let us now substitute (\ref{Es2}) for (\ref{Es1}); then
%
\begin{eqnarray}\label{EE2}
\dif X_{t,s}(x)=-[\mE\varphi(X_{s,0}(y))]_{y=X_{t,s}(x)}\,\dif s+\sqrt
{2\nu}\,\dif W_s,
\nonumber
\\[-8pt]
\\[-8pt]
\eqntext{s\in[t,0],\  X_{t,t}(x)=x.}
\end{eqnarray}
As the Markov property holds, one can write the above equation as a
closed form,
%
\begin{equation}\label{Es4}
\dif X_{t,s}(x)=-\mE^{\sF_{t,s}}\varphi(X_{t,0}(x))\,\dif s+\sqrt
{2\nu
}\,\dif W_s,\qquad s\in[t,0],\
  X_{t,t}(x)=x,\hspace*{-35pt}
\end{equation}
where $\sF_{t,s}:=\sigma\{W_r-W_t\dvtx  r\in[t,s]\}$, and $\mE^{\sF_{t,s}}$
denotes the conditional expectation
with respect to $\sF_{t,s}$. The question is this: Suppose that the
stochastic equation~(\ref{Es4}) admits
a unique solution family $\{X_{t,s}(x)\dvtx  t\leq s\leq0, x\in\mR^d\}$.
Does $u_t(x)$
defined by (\ref{Es2}) solve Burgers' equation (\ref{Es3})? To answer
this question, the key point is to
establish the following Markov property: for all $t_1\leq t_2\leq
t_3\leq0$ and $x\in\mR^d$,
%
\begin{equation}\label{Mar}
\mE^{\sF_{t_1,t_2}}(\varphi(X_{t_1,t_3}(x)))=\mE(\varphi
(X_{t_2,t_3}(y)))|_{y=X_{t_1,t_2}(x)}\qquad   \mbox{a.s.}
\end{equation}
so that equation (\ref{Es4}) can be written back to (\ref{EE2}). This
is not obvious
since SDE (\ref{Es4}) involves a conditional expectation operator.
On the other hand, one can replace the Brownian motion in equation
(\ref{Es4})
by an $\alpha$-stable process, as is done in~\cite{Zh2}, so that we can
give a probabilistic explanation
for the Burgers equation (\ref{Es3}).

Basing on this simple observation, in this paper we are mainly
concerned about
the following general stochastic functional differential equation
(abbreviated as SFDE) driven by a L\'evy process $(L_t)_{t\leq0}$:
%
\begin{eqnarray}\label{SDE0}
\dif X_{t,s}(x)=G_s(X_{t,s-}(x),\mE^{\sF_{s-}} (\phi_s(X_{t,\cdot
}(x))))\,\dif L_s,
\nonumber
\\[-8pt]
\\[-8pt]
\eqntext{ s\in[t,0],\ X_{t,t}(x)=x,}
\end{eqnarray}
where $\sF_s:=\sigma\{L_{s'}-L_{s''}\dvtx  s''<s'\leq s\}$, $G$ and $\phi$
are some Lipschitz functionals (see below).
In Section~\ref{sec2}, we are devoted to proving the existence and uniqueness of
a short time solution
as well as the Markov property (\ref{Mar}) for equation (\ref{SDE0})
under Lipschitz assumptions
on $G$ and $\phi$. Moreover, a locally maximal solution is also achieved.
Since L\'evy processes usually have poor integrability, we have to
carefully treat the
big jump part of L\'evy processes. Compared with the classical argument
in Freidlin~\cite{Fr}, it seems that SFDE (\ref{SDE0}) is easier to
handle since it is a closed equation.

Next, in Section~\ref{sec3} we apply our results to a class of quasi-linear
partial integro-differential equations
(abbreviated as PIDE) and
obtain the existence of short time solutions. Here, we discuss two
cases: $G$ and $\phi$ satisfy
linear growth conditions, but L\'evy processes have finite moments of
arbitrary orders; $G$ and $\phi$ are bounded,
but equation (\ref{SDE0}) has a constant coefficient in the big jump part.
This is natural since only the big jump is related to
the moment of L\'evy processes.

In Section~\ref{sec4}, we turn to the investigation of the following system of
semi-linear PIDEs
(nonlinear transport equation):
%
\begin{equation}\label{Tr}
\cases{
\partial _t u_t+\cL_0 u_t+\bigl(G_t(x, u_t)\cdot\nabla\bigr)
u_t+F_t(x,u_t)=0,\vspace*{2pt}\cr
(t,x)\in\mR_-\times\mR^d, \qquad  u_0(x)=\varphi(x)\in\mR^m,}
\end{equation}
where $\cL_0$ is the generator of the L\'evy process given by (\ref
{Ge}) below. It is observed that
the following scalar conservation law equation can be written as the
above form:
%
\begin{equation}\label{Sc}
\cases{
\partial _t u_t+\cL_0 u_t+\div(g_t(x,
u_t))+f_t(x,u_t)=0,\vspace*{2pt}\cr
(t,x)\in\mR_-\times\mR^d,  \qquad u_0(x)=\varphi(x)\in\mR.}
\end{equation}
In particular, the one-dimensional fractal Burgers equation (\ref{Es3})
takes the above form.
In equation (\ref{Tr}), since there are not any analytic properties to
be imposed on $\cL_0$, one can not appeal to the Duhamel formula or the
energy method to give
an analytic treatment. In this situation, the probabilistic approach
seems to be quite suitable.
In fact, by using purely probabilistic argument, we shall prove in
Theorem~\ref{Main} below that
PIDE (\ref{Tr}) admits a unique maximal weak solution in the class of
bounded Lipschitz functions.
In the nondegenerate case (corresponding to the subcritical case for
$\cL_0=\Delta^{\alpha/2}$ with $\alpha\in(1,2]$),
the existence of global solutions is also obtained by applying some
gradient estimates.
We mention that for the one-dimensional Burgers equation (\ref{EL5}),
it has been proved in~\cite{Ki-Na-Sh} that the global analytic solution
does exist for $\alpha\in[1,2]$,
and the finite time blow up solution also exists for $\alpha\in
(0,1)$. However,
in the critical case of $\alpha=1$, the existence of global solutions
for the general equation (\ref{Tr}) is left open.

We conclude this introduction by introducing the following conventions:
The letter $C$ with or without subscripts
will denote a positive constant, whose value may change in different
places. If we write $T=T(K_1,K_2,\ldots)$,
this means that $T$ depends only on these indicated arguments.

\section{A stochastic functional differential equation: Short time existence}\label{sec2}

\subsection{General facts about L\'evy processes}\label{sec2.1}

Let $(L_t)_{t\in\mR}$ be an $\mR^m$-valued L\'evy process on the real
line and defined on
some complete probability space $(\Omega,\sF,P)$, which means that:
\begin{itemize}
\item $(L_t)_{t\in\mR}$ has independent and stationary increments,
that is, for all
$-\infty<t_1<t_2<\cdots<t_n<+\infty$,
the random variables $(L_{t_2}-L_{t_1},\ldots, L_{t_n}-L_{t_{n-1}})$
are independent, and the distribution of
$L_{t+s}-L_s$ does not depend on $s$.

\item For $P$-almost all $\omega\in\Omega$, the mapping
$t\mapsto
L_t(\omega)$
is right-continuous and has left-limit (also called c\`adl\`ag in French).
\end{itemize}

Let $\sN$ be the total of all $P$-null sets. For $-\infty\leq
t<s<+\infty$, define
\[
\sF_{t,s}:=\sigma\{L_r-L_{r'}; r,r'\in(t,s]\}\vee\sN.
\]
By the independence of increments of the L\'evy process, it is easy to
see that
for $-\infty\leq t_1<t_2<t_3<+\infty$, $\sF_{t_1,t_2}$ and $\sF
_{t_2,t_3}$ are independent.
For simplicity of notation, we write
\[
\sF_s=\sF_{-\infty,s}, \qquad \sF_{s-}:=\bigvee_{t<s}\sF_t.
\]
It is clear that $\sF_t\subset\sF_s$ if $t<s$, and $s\mapsto\sF_{s-}$
is left-continuous. Moreover,
$L_{s-}$ is $\sF_{s-}$-measurable.
Throughout this paper, we shall work on the negative time axes $\mR
_-:=(-\infty,0]$.

\begin{remark}\label{Re1}
For any measurable process $\eta_s\in L^1(\Omega,\sF_0,P)$, $s\leq0$,
by the predictable projection theorem
(cf.~\cite{Re-Yo}, page~173, Theorem~5.3), there always exists a
predictable version of
$s\to\mE(\eta_s|\sF_{s-})$, which will be denoted by $\mE^{\sF
_{s-}}(\eta_s)$.
Moreover, for any $\xi\in L^1(\Omega,\sF_0,P)$, by the regularization
theorem of martingales
(cf.~\cite{Re-Yo}, page~64, Proposition~2.7 and page~65,\vadjust{\goodbreak} Theorem~2.9),
we have
\[
\lim_{s\uparrow t}\mE^{\sF_{s-}}(\xi)=\mE^{\sF_{t-}}(\xi)=\mE
^{\sF
_t}(\xi)\qquad \mbox{a.s.},
\]
where the second equality follows by $P\{L_t=L_{t-}\}=1$ and a monotone
class argument.
\end{remark}

By the L\'evy--Khintchine formula (cf.~\cite{Ap}, page~109, Corollary~2.4.20),
the characteristic function of $L_t$ is given by
%
\begin{eqnarray}\label{Lp9}
\quad\mE(e^{i\xi\cdot L_t})
&=&\exp\biggl\{t\biggl[ib\cdot\xi-\xi^{\mathrm{t}} A\xi+\int_{\mR
^m}\bigl[e^{i\xi
\cdot z}-1
-i\xi\cdot z 1_{|z|\leq1}\bigr]\nu(\dif z)\biggr]\biggr\}
\nonumber
\\[-8pt]
\\[-8pt]
\nonumber
&=:&e^{t\Psi(\xi)},
\end{eqnarray}
where $\Psi(\xi)$ is a complex-valued function called the symbol of
$(L_t)_{t\leq0}$, and
$b\in\mR^m$, $A\in\mR^m\times\mR^m$ is a positive definite and
symmetric matrix,
$\nu$ is a L\'evy measure on $\mR^m$, that is, $\nu\{0\}=0$ and
%
\begin{equation}\label{Levy}
\int_{\mR^m}1\wedge|z|^2\nu(\dif z)<+\infty.
\end{equation}
We call
%
\begin{equation}\label{Ch}
\sA:=(b,A,\nu)
\end{equation}
the characteristic triple of $L_t$. If $b=0, A=0$ and $\nu(\dif
z)=\frac
{\dif z}{|z|^{m+\alpha}}$, where $\alpha\in(0,2)$,
then $L_t$ is the $\alpha$-stable process with the L\'evy exponent
$c_{m,\alpha}|\xi|^\alpha$,
and its generator is the fractional Laplacian $\Delta^{\alpha/2}$ by
multiplying a constant~$c_{m,\alpha}'$.

By the L\'evy--It\^o decomposition (cf.~\cite{Ap}, page~108, Theorem~2.4.16), $L_t$ can be written as
%
\begin{equation}\label{Le-Ito}
L_t=bt+W^A_t+\int_{|z|\leq1}z\tilde N(t,\dif z)+\int
_{|z|>1}zN(t,\dif
z),
\end{equation}
where $W^A_t$ is a Brownian motion with covariance matrix $A=(a_{ij})$,
$N(t,\dif z)$ is the Poisson random point measure associated with
$(L_t)_{t\leq0}$ given by
\[
N(t,\Gamma):=\sum_{t<s\leq0}1_{\Gamma}(L_s-L_{s-}),\qquad \Gamma\in\cB
(\mR^m)
\]
and $\tilde N(t,\dif z):=N(t,\dif z)-t\nu(\dif z)$ is the compensated
random martingale measure.
Here, $(W^A_t)_{t\leq0}$ and $(N(t,\dif z))_{t\leq0}$ are
independent. The generator of $L_t$ is given by
%
\begin{eqnarray}\label{Ge}
\cL_0 u(x)&=&\frac{1}{2}a_{ij}\partial _i\partial _j u+b_i\partial
_i u
\nonumber
\\[-8pt]
\\[-8pt]
\nonumber
&&{}+\int_{\mR
^m}
\bigl[u(x+z)-u(x)-1_{|z|<1}
\partial _i u(x)z_i\bigr]\nu(\dif z).
\end{eqnarray}
Here and after, we use the usual convention for summation: the same
index in a product will be summed
automatically.\vadjust{\goodbreak}

In the following, we denote by $\mD$ the space of all c\`adl\`ag
functions from $\mR_-$ to~$\mR^d$, which is
endowed with the locally uniform metric $\rho$. Notice that this metric
is complete, but not separable.
For given $t<0$ and c\`adl\`ag function $f:[t,0]\to\mR^d$, we extend
$f$ to $\mR_-$ in a natural manner by putting
$f(s)=f(t)$ for $s<t$ so that $f\in\mD$.

\subsection{A general case}\label{sec2.2}

In this subsection, we consider the following general SFDE in $\mR^d$ driven
by the L\'evy process $(L_s)_{s\leq0}$:
%
\begin{equation}\label{SDE}
X_{t,s}=\xi+\int_{(t,s]}G_r(X_{t,r-},\mE^{\sF_{r-}} (\phi
_r(X_{t,\cdot
})))\,\dif L_r,\qquad  t\leq s\leq0,
\end{equation}
where $\xi\in\sF_t$, $G\dvtx \mR_-\times\mR^d\times\mR^k\to\mR
^d\times\mR^m$
is a measurable function,
and $\phi\dvtx \mR_-\times\mD\to\mR^k$ is a uniformly Lipschitz continuous
functional in the sense that
%
\begin{equation}\label{EL1}
\|\phi\|_{\mathrm{Lip}}:=
\sup_{s\in\mR_-}\sup_{\omega\not=\omega'\in\mD}
\frac{|\phi_s(\omega)-\phi_s(\omega')|}{\rho(\omega,\omega
')}<+\infty
,
\end{equation}
where $\rho(\omega,\omega'):=\sum_{n}2^{-n}(1\wedge\sup
_{s\in
[-n,0]}|\omega(s)-\omega'(s)|)$ is the locally uniform metric on
$\mD$.

The definition about the solutions to equation (\ref{SDE}) is given as follows:
%
\begin{definition}\label{Def}
For fixed $t<0$ and $\xi\in\sF_t$, an $(\sF_s)$-adapted c\`adl\`ag
stochastic process
$X_s=:X_{t,s}(\xi)$ is called a solution of equation (\ref{SDE}) if for
all $s\in[t,0]$,
\[
X_s=\xi+\int_{(t,s]}G_r(X_{r-},\mE^{\sF_{r-}}(\phi_r(X_\cdot
)))\,\dif
L_r\qquad \mbox{a.s.}
\]
For $T<0$, we say that equation (\ref{SDE}) is (uniquely) solvable on $(T,0]$
(or $[T,0]$) if for all $t\in(T,0]$ (or $t\in[T,0]$)
and $\xi\in\sF_t$, equation (\ref{SDE}) has a (unique) solution
starting from $\xi$ at time $t$.
\end{definition}

\begin{remark}
In this definition, it has been assumed that $\phi_r(X_\cdot)\in
L^1(\Omega, \sF_0,P)$
so that $\mE^{\sF_{r-}}(\phi_r(X_\cdot))$
makes sense by Remark~\ref{Re1}, and further the stochastic integral
with respect to the L\'evy process
in the definition makes sense.
\end{remark}

Below, we make the following assumptions on the coefficients and the L\'
evy measure $\nu$:
\begin{longlist}[(\textbf{H}$^\beta_\nu$)]
\item[(\textbf{H}$_G$)] For some $K_0, K_1>0$ and all $s\leq0$, $x,x'\in\mR^d$,
$u,u'\in
\mR^k$,
\[
|G_s(0,0)|\leq K_0, \qquad |G_s(x,u)-G_s(x',u')|\leq K_1(|x-x'|+|u-u'|).
\]
\item[(\textbf{H}$^\beta_\nu$)] For some $\beta>0$,
\[
\int_{|z|\geq1}|z|^\beta\nu(\dif z)<+\infty.\vadjust{\goodbreak}
\]
\end{longlist}

\begin{remark}
Condition (\textbf{H}$^\beta_\nu$), which is a restriction on the big jump
of the L\'evy process,
is equivalent to saying that the $\beta$-order moment of the L\'evy
process is finite; cf.~\cite{Sa}, Theorem~25.3. It should be noticed that
for $\alpha$-stable process, condition (\textbf{H}$^\beta_\nu$) is
satisfied only for any $\beta<\alpha$.
\end{remark}

Now we prove the following result about the existence and uniqueness of
solutions for equation (\ref{SDE})
in a short time.
%
\begin{theorem}\label{Th1}
Assume that \textup{(\textbf{H}$_G$)} and \textup{(\textbf{H}$^\beta_\nu$)} hold for some
$\beta
>1$, and $\phi$ is a
Lipschitz continuous functional on $\mD$; see (\ref{EL1}).
Then there exists a time $T=T(K_1,\sA,\beta,\|\phi\|_{\mathrm{Lip}})<0$ such that
equation (\ref{SDE}) is uniquely solvable on $[T,0]$ for any $L^\beta
$-integrable initial value $\xi\in\sF_t$
in the sense of Definition~\ref{Def}, and for some $C=C(T,K_0)$ and any
$t\in[T,0]$,
%
\begin{equation}\label{EL6}
\mE\Bigl(\sup_{s\in[t,0]}|X_{t,s}(\xi)|^\beta\Bigr)\leq C\mE
|\xi|^\beta
.
\end{equation}
Moreover, if $\xi=x\in\mR^d$ is nonrandom, then for any $t\in[T,0)$,
the unique solution $X_{t,s}$ is $\sF_{t,s}$-measurable for all $s\in[t,0]$.
\end{theorem}

\begin{pf}
We prove the theorem for $\beta\in(1,2)$. For $\beta\geq2$, the proof
is similar and simpler.
Fix $t<0$, which will be determined below.
For $\xi\in L^\beta(\Omega,\sF_t,P)$, set $X^{(0)}_{t,s}\equiv\xi
$, and
let $X^{(n)}_{t,s}$ be the Picard iteration sequence defined as
follows: for $n\in\mN$,
%
\begin{equation}\label{Ep3}
X^{(n)}_{t,s}=\xi+\int_{(t,s]}G_r\bigl(X^{(n-1)}_{t,r-},\mE^{\sF_{r-}}
\bigl(\phi_r\bigl(X^{(n-1)}_{t,\cdot}\bigr)\bigr)\bigr)\,\dif L_r.
\end{equation}
Set
\[
Z^{(n)}_{t,s}:=X^{(n+1)}_{t,s}-X^{(n)}_{t,s}.
\]
Using the L\'evy--It\^o decomposition (\ref{Le-Ito}), one can write
\begin{eqnarray*}
Z^{(n)}_{t,s}&=&\int_{(t,s]} \int_{|z|<1}\cG^{(n)}_r\cdot z\tilde
N(\dif r,\dif z)\\
&&{}+\int_{(t,s]} \int_{|z|\geq1}\cG^{(n)}_r\cdot z N(\dif r,\dif z)\\
&&{} +\int_{(t,s]}\cG^{(n)}_r\cdot b\,\dif r+\int_{(t,s]}\cG
^{(n)}_r\,\dif
W^A_r\\
&=:&I^{(n)}_1(s)+I^{(n)}_2(s)+I^{(n)}_3(s)+I^{(n)}_4(s),
\end{eqnarray*}
where
\[
\cG^{(n)}_r:=G_r\bigl(X^{(n)}_{t,r-},\mE^{\sF_{r-}}
\bigl(\phi_r\bigl(X^{(n)}_{t,\cdot}\bigr)\bigr)\bigr)-G_r\bigl(X^{(n-1)}_{t,r-},\mE^{\sF_{r-}}
\bigl(\phi_r\bigl(X^{(n-1)}_{t,\cdot}\bigr)\bigr)\bigr).
\]
By Burkholder's inequality (cf.~\cite{Ka}, Theorem~23.12) and Young's
inequality, thanks to $\beta\in(1,2)$,
we have that for any $\varepsilon\in(0,1)$,
%
\begin{eqnarray*}
&&\mE\Bigl(\sup_{r\in[t,0]}\bigl|I^{(n)}_1(r)\bigr|^\beta\Bigr)\\
&&\qquad\leq C
\mE\biggl(\int_{(t,0]} \int_{|z|<1}\bigl|\cG^{(n)}_r\cdot z\bigr|^2N(\dif
r,\dif
z)\biggr)^{\beta/2} \\
&&\qquad\leq C\mE\biggl(\sup_{r\in[t,0]}\bigl|\cG^{(n)}_r\bigr|^{2-\beta}\int
_{(t,0]}
\int_{|z|<1}\bigl|\cG^{(n)}_r\bigr|^\beta
\cdot|z|^2N(\dif r,\dif z)\biggr)^{\beta/2} \\
&&\qquad\leq\varepsilon\mE\Bigl(\sup_{r\in[t,0]}\bigl|\cG^{(n)}_r\bigr|^\beta\Bigr)
+C_\varepsilon\mE\biggl(\int_{(t,0]} \int_{|z|<1}\bigl|\cG^{(n)}_r\bigr|^\beta
\cdot|z|^2N(\dif r,\dif z)\biggr) \\
&&\qquad=\varepsilon\mE\Bigl(\sup_{r\in[t,0]}\bigl|\cG^{(n)}_r\bigr|^\beta\Bigr)
+C_\varepsilon\mE\biggl(\int_{[t,0]} \int_{|z|<1}\bigl|\cG^{(n)}_r\bigr|^\beta
\cdot|z|^2\nu(\dif z)\,\dif r\biggr) \\
&&\qquad\leq\biggl(\varepsilon+C_\varepsilon|t|\int_{|z|<1}|z|^2\nu(\dif
z)\biggr)
\mE\Bigl(\sup_{r\in[t,0]}\bigl|\cG^{(n)}_r\bigr|^\beta\Bigr).
\end{eqnarray*}
Here and below, the constant $C$ or $C_\varepsilon$ is independent of $t$
and $n$.
For $I^{(n)}_2(s)$, by It\^o's formula, we have
\begin{eqnarray*}
&&\mE\Bigl(\sup_{r\in[t,s]}\bigl|I^{(n)}_2(r)\bigr|^\beta\Bigr)\\
&&\qquad\leq\mE
\biggl(\int
_{(t,s]} \int_{|z|\geq1}\bigl|\bigl|I^{(n)}_2(r-)
+\cG^{(n)}_r\cdot z\bigr|^\beta-\bigl|I^{(n)}_2(r-)\bigr|^\beta\bigr|N(\dif r,\dif
z)\biggr)\\
&&\qquad=\mE\biggl(\int_{(t,s]} \int_{|z|\geq1}\bigl|\bigl|I^{(n)}_2(r-)
+\cG^{(n)}_r\cdot z\bigr|^\beta-\bigl|I^{(n)}_2(r-)\bigr|^\beta\bigr|\nu(\dif
z)\,\dif
r\biggr)\\
&&\qquad\leq C \mE\biggl(\int_{(t,s]}\bigl|I^{(n)}_2(r)\bigr|^\beta\,\dif r
\biggr)+C
\biggl(\int_{|z|\geq1}|z|^\beta\nu(\dif z)\biggr)
\mE\biggl(\int_{(t,0]}\bigl|\cG^{(n)}_r\bigr|^\beta\,\dif r\biggr),
\end{eqnarray*}
which then implies that by (\textbf{H}$^\beta_\nu$) and Gronwall's inequality,
\[
\mE\Bigl(\sup_{r\in[t,0]}\bigl|I^{(n)}_2(r)\bigr|^\beta\Bigr)\leq
C|t|\mE\Bigl(\sup_{r\in[t,0]}\bigl|\cG^{(n)}_r\bigr|^\beta\Bigr).
\]
Similarly, we have
\[
\mE\Bigl(\sup_{r\in[t,0]} \bigl|I^{(n)}_3(r)\bigr|^\beta\Bigr)\leq
(|t|\cdot|b|)^\beta\mE\Bigl(\sup_{r\in[t,0]}\bigl|\cG^{(n)}_r\bigr|^\beta
\Bigr),
\]
and for any $\varepsilon\in(0,1)$,
\[
\mE\Bigl(\sup_{r\in[t,0]}\bigl|I^{(n)}_4(r)\bigr|^\beta\Bigr)\leq
(\varepsilon
+C_\varepsilon|t|)
\mE\Bigl(\sup_{r\in[t,0]}\bigl|\cG^{(n)}_r\bigr|^\beta\Bigr).
\]
Combining the above calculations, we obtain that for any $\varepsilon\in(0,1)$,
%
\begin{equation}\label{EP5}
\mE\Bigl(\sup_{r\in[t,0]}\bigl|Z^{(n)}_{t,r}\bigr|^\beta\Bigr)\leq
(\varepsilon
+C_\varepsilon|t|)\cdot
\mE\Bigl(\sup_{r\in[t,0]}\bigl|\cG^{(n)}_r\bigr|^\beta\Bigr).
\end{equation}
Noticing that by (\textbf{H}$_G$),
\[
\bigl|\cG^{(n)}_r\bigr|\leq K_1\Bigl(\bigl|Z^{(n-1)}_{t,r-}\bigr|+\|\phi\|_{\mathrm{Lip}}\mE
^{\sF_{r-}}
\Bigl(\sup_{s\in[t,0]}\bigl|Z^{(n-1)}_{t,s}\bigr|\Bigr)\Bigr),
\]
and in view of $\beta>1$, we further have by Doob's maximal inequality,
\[
\mE\Bigl(\sup_{s\in[t,0]}\bigl|Z^{(n)}_{t,s}\bigr|^\beta\Bigr)
\leq(\varepsilon+C_\varepsilon|t|)C_0\mE\Bigl(\sup_{s\in
[t,0]}\bigl|Z^{(n-1)}_{t,s}\bigr|^\beta\Bigr).
\]
Now, let us choose
\[
\varepsilon=\frac{1}{4C_0}\quad\mbox{and}\quad
T:=-\frac{1}{4C_\varepsilon C_0},
\]
and then for all $t\in[T,0]$,
%
\begin{equation}\label{Eo1}
\qquad\mE\Bigl(\sup_{s\in[t,0]}\bigl|Z^{(n)}_{t,s}\bigr|^\beta\Bigr)\leq
\frac{1}{2}\mE\Bigl(\sup_{s\in[t,0]}\bigl|Z^{(n-1)}_{t,s}\bigr|^\beta
\Bigr)\leq
\cdots\leq\frac{1}{2^n}
\mE\Bigl(\sup_{s\in[t,0]}\bigl|Z^{(0)}_{t,s}\bigr|^\beta\Bigr).
\end{equation}
On the other hand, notice that
\[
Z^{(0)}_{t,s}=X^{(1)}_{t,s}-\xi=\int_{(t,s]}G_r\bigl(X^{(1)}_{t,r-},\mE
^{\sF
_{r-}}(\phi_r(\xi))\bigr)\,\dif L_r.
\]
As above, and using Gronwall's inequality, it is easy to derive that
%
\begin{equation}\label{EL7}
\mE\Bigl(\sup_{s\in[t,0]}\bigl|Z^{(0)}_{t,s}\bigr|^\beta\Bigr)\leq C\mE
|\xi
|^\beta.
\end{equation}
Hence, there exists an ($\sF_s$)-adapted and c\`adl\`ag stochastic
process $X_{t,s}$ such that
%
\begin{equation}\label{EL8}
\lim_{n\to\infty}\mE\Bigl(\sup_{s\in
[t,0]}\bigl|X^{(n)}_{t,s}-X_{t,s}\bigr|^\beta
\Bigr)=0.
\end{equation}
By taking limits for equation (\ref{Ep3}), it is easy to see that
$X_{t,s}$ solves SFDE (\ref{SDE}).
Moreover, estimate (\ref{EL6}) follows from (\ref{Eo1}), (\ref{EL7})
and (\ref{EL8}).
The uniqueness is clear from the above proof.

Suppose now that $\xi=x$ is nonrandom. From the Picard iteration
(\ref
{Ep3}), one sees that for each $n\in\mN$
and $s\in(t,0]$, $X^{(n)}_{t,s}$ is $\sF_{t,s}$-measurable. Indeed,
suppose that
$X^{(n-1)}_{t,s}$ is $\sF_{t,s}$-measurable for each $s\in(t,0]$, and
then it is clear that $\phi_r(X^{(n-1)}_{t,\cdot})$
is independent of $\sF_t$. Noticing that for $r>t$,
$\sF_{r-}=\sF_{t,r-}\vee\sF_t$ and $\sF_{t,r-}$ is independent of
$\sF
_t$, we have
\[
\mE^{\sF_{r-}} \bigl(\phi_r\bigl(X^{(n-1)}_{t,\cdot}\bigr)\bigr)=\mE^{\sF_{t,r-}}
\bigl(\phi
_r\bigl(X^{(n-1)}_{t,\cdot}\bigr)\bigr).
\]
By induction method, starting from equation (\ref{Ep3}) with $\xi=x$,
one finds that
$X^{(n)}_{t,s}$ is also $\sF_{t,s}$-measurable for each $s\in(t,0]$.
So, the limit $X_{t,s}$ is also $\sF_{t,s}$-measurable.
\end{pf}
%
\begin{remark}
In this theorem, if $G_s(x,u)=G_s(x)$ does not depend on~$u$,
then the short time solution can be extended to any large time by the
usual time shift technique.
\end{remark}

\subsection{A special case}\label{sec2.3}

In Theorem~\ref{Th1}, since we require $\beta>1$, the result rules out
the $\alpha$-stable process
with $\alpha\in(0,1]$. In this subsection, we drop assumption (\textbf{H}$^\beta_\nu$) in Theorem~\ref{Th1},
and consider the following special form:
%
\begin{eqnarray}\label{SDE2}
X_{t,s}&=&\xi+\int_{(t,s]} \int_{|z|<1}G_r(X_{t,r-},\mE^{\sF
_{r-}}(\phi
_r(X_{t,\cdot})))\cdot z\tilde N(\dif r,\dif z)\nonumber\\[-2pt]
&&{}+\int_{(t,s]} \int_{|z|\geq1}z N(\dif r,\dif z)
\nonumber
\\[-9pt]
\\[-9pt]
\nonumber
&&{} +\int_{(t,s]}G_r(X_{t,r-},\mE^{\sF_{r-}}(\phi_r(X_{t,\cdot
})))\cdot b\,\dif r
\\[-2pt]
&&{}+\int_{(t,s]}G_r(X_{t,r-},\mE^{\sF_{r-}}(\phi_r(X_{t,\cdot})))\,\dif
W^A_r,\nonumber
\end{eqnarray}
where $\xi\in\sF_t$. In this equation, the big jump part has a constant
coefficient. In order to make
sense for the integrals, we need to assume that $G$ and $\phi$ are
bounded. We have:
%
\begin{theorem}\label{Th3}
In addition to \textup{(\textbf{H}$_G$)}, we assume that $G$ is bounded, and~$\phi$
is a bounded
Lipschitz continuous functional on $\mD$. Then there exists a time
$T=T(K_1,\sA,\|\phi\|_{\mathrm{Lip}})<0$ such that
SFDE (\ref{SDE2}) is uniquely solvable on $[T,0]$.
Moreover, if $\xi=x\in\mR^d$ is nonrandom, then for any $t\in[T,0)$,
the unique solution $X_{t,s}$ is $\sF_{t,s}$-measurable for all $s\in[t,0]$.
\end{theorem}

\begin{pf}
For $t<0$ and $\xi\in\sF_t$, set $X^{(0)}_{t,s}\equiv\xi$, and
let $X^{(n)}_{t,s}$ be the Picard iteration sequence defined as follows:
%
\begin{eqnarray}\label{Ep33}
X^{(n)}_{t,s}&=&\xi+\int_{(t,s]} \int_{|z|\geq1}z N(\dif r,\dif z)
\nonumber\\[-2pt]
&&{}+\int_{(t,s]} \int_{|z|<1}G_r\bigl(X^{(n-1)}_{t,r-},\mE^{\sF_{r-}}\bigl(\phi
_r\bigl(X^{(n-1)}_{t,\cdot}\bigr)\bigr)\bigr)\cdot z\tilde N(\dif r,\dif z)
\nonumber
\\[-9pt]
\\[-9pt]
\nonumber
&&{} +\int_{(t,s]}G_r\bigl(X^{(n-1)}_{t,r-},\mE^{\sF_{r-}}\bigl(\phi
_r\bigl(X^{(n-1)}_{t,\cdot}\bigr)\bigr)\bigr)\cdot b\,\dif r
\\[-2pt]
&&{}+\int_{(t,s]}G_r\bigl(X^{(n-1)}_{t,r-},\mE^{\sF_{r-}}\bigl(\phi
_r\bigl(X^{(n-1)}_{t,\cdot}\bigr)\bigr)\bigr)\,\dif W^A_r.\nonumber
\end{eqnarray}
Set
\[
Z^{(n)}_{t,s}:=X^{(n+1)}_{t,s}-X^{(n)}_{t,s}.
\]
Then
\[
Z^{(n)}_{t,s}=\int_{(t,s]} \int_{|z|<1}\cG^{(n)}_r\cdot z\tilde
N(\dif
r,\dif z)
+\int_{(t,s]}\cG^{(n)}_r\cdot b\,\dif r+\int_{(t,s]}\cG^{(n)}_r\,\dif W^A_r,
\]
where
\[
\cG^{(n)}_r:=G_r\bigl(X^{(n)}_{t,r-},\mE^{\sF_{r-}}
\bigl(\phi_r\bigl(X^{(n)}_{t,\cdot}\bigr)\bigr)\bigr)-G_r\bigl(X^{(n-1)}_{t,r-},\mE^{\sF
_{r-}}\bigl(\phi
_r\bigl(X^{(n-1)}_{t,\cdot}\bigr)\bigr)\bigr).
\]
Notice that $|\cG^{(n)}_r|\leq2\|G\|_\infty$, and by (\textbf{H}$_G$),
\[
\bigl|\cG^{(n)}_r\bigr|^2\leq2K^2_1\Bigl(\bigl|Z^{(n-1)}_{t,r-}\bigr|^2+\|\phi\|
^2_{\mathrm{Lip}}\Bigl(\mE^{\sF_{r-}}
\Bigl(\sup_{s\in[t,0]}\bigl|Z^{(n-1)}_{t,s}\bigr|\Bigr)\Bigr)^2
\Bigr)=:\Phi^{(n)}_r.
\]
By Burkholder's inequality and (\ref{Levy}), we have
\begin{eqnarray*}
\mE\biggl(\sup_{s\in[t,0]}\biggl|\int_{(t,s]} \int_{|z|<1}\cG
^{(n)}_r\cdot z\tilde N(\dif r,\dif z)\biggr|^2\biggr)
&\leq &C\mE\biggl(\int_{(t,0]} \int_{|z|<1}\bigl|\cG^{(n)}_r\cdot
z\bigr|^2N(\dif
r,\dif z)\biggr)
\\
&\leq &C|t|\mE\Bigl(\sup_{r\in[t,0]}\Phi^{(n)}_r\Bigr).
\end{eqnarray*}
Here and below, the constant $C$ is independent of $t$ and $n$.
Similarly, we have
\[
\mE\biggl(\sup_{s\in[t,0]}\biggl|\int_{(t,s]}\cG^{(n)}_r\cdot
b\,\dif
r\biggr|^2\biggr)\leq
C|t|^2\mE\Bigl(\sup_{r\in[t,0]}\Phi^{(n)}_r\Bigr)
\]
and
\[
\mE\biggl(\sup_{s\in[t,0]}\biggl|\int_{(t,s]}\cG^{(n)}_r\,\dif
W^A_r
\biggr|^2\biggr)\leq
C|t|\mE\Bigl(\sup_{r\in[t,0]}\Phi^{(n)}_r\Bigr).
\]
Combining the above calculations and by Doob's maximal inequality, we obtain
\[
\mE\Bigl(\sup_{s\in[t,0]}\bigl|Z^{(n)}_{t,s}\bigr|^2\Bigr)\leq C_0|t|\mE
\biggl(\sup_{r\in[t,0]}\Phi^{(n)}_r\biggr)
\leq C_1|t|\mE\Bigl(\sup_{s\in[t,0]}\bigl|Z^{(n-1)}_{t,s}\bigr|^2\Bigr).
\]
Now, let us choose
\[
T:=-\frac{1}{2C_1},
\]
and then for all $t\in[T,0]$,
\[
\mE\Bigl(\sup_{r\in[t,0]}\bigl|Z^{(n)}_{t,r}\bigr|^2\Bigr)\leq
\frac{1}{2}\mE\Bigl(\sup_{r\in[t,0]}\bigl|Z^{(n-1)}_{t,r}\bigr|^2
\Bigr)\leq\cdots
\leq
\frac{1}{2^n}\mE\Bigl(\sup_{r\in[t,0]}\bigl|Z^{(1)}_{t,r}\bigr|^2
\Bigr)\leq\frac{C}{2^n}.
\]
Hence, there exists an ($\sF_s$)-adapted and c\`adl\`ag stochastic
process $X_{t,s}$ such that
\[
\lim_{n\to\infty}\mE\Bigl(\sup_{s\in
[t,0]}\bigl|X^{(n)}_{t,s}-X_{t,s}\bigr|^2\Bigr)=0.
\]
By taking limits for equation (\ref{Ep33}), it is easy to see that
$X_{t,s}$ solves SFDE (\ref{SDE2}).
The remaining proof is the same as in Theorem~\ref{Th1}.
\end{pf}

\subsection{Markov property}\label{sec2.4}
In this subsection, we prove the Markov property for the solutions of
equations (\ref{SDE})
and (\ref{SDE2}), which is crucial for the development of the next section.

We first show the continuous dependence of the solutions with respect
to the initial values.
%
\begin{proposition}\label{Pro1}
In the situation of Theorem~\ref{Th1}, for $t\in[T,0]$,
let $\xi^{(n)},\xi\in L^\beta(\Omega,\sF_t,P)$. If $\xi^{(n)}$
converges to $\xi$ in probability as $n\to\infty$, then
$X^{(n)}_{t,s}$ converges to $X_{t,s}$ uniformly with respect to $s\in
[t,0]$ in probability
as $n\to\infty$, where $\{X^{(n)}_{t,s}; t\leq s\leq0\}$ and $\{
X_{t,s}; t\leq s\leq0\}$
are the solutions of SFDE (\ref{SDE}) corresponding to the initial
values $\xi^{(n)}$ and $\xi$.
\end{proposition}

\begin{pf}
Define
\[
A_n:=\bigl\{\bigl|\xi^{(n)}-\xi\bigr|\leq1\bigr\}\in\sF_t.
\]
Then we can write
\begin{eqnarray*}
1_{A_n}\bigl(X^{(n)}_{t,s}-X_{t,s}\bigr)&=&1_{A_n}\bigl(\xi^{(n)}-\xi\bigr)+
\int_{(t,s]} \int_{|z|<1}1_{A_n}\cG^{(n)}_r\cdot z\tilde N(\dif
r,\dif
z)\\
&&{} +\int_{(t,s]} \int_{|z|\geq1}1_{A_n}\cG^{(n)}_r\cdot z
N(\dif
r,\dif z)\\
&&{} +\int_{(t,s]}1_{A_n}\cG^{(n)}_r\cdot b\,\dif r+\int
_{(t,s]}1_{A_n}\cG^{(n)}_r\,\dif W^A_r,
\end{eqnarray*}
where
\[
\cG^{(n)}_r:=G_r\bigl(X^{(n)}_{t,r-},\mE^{\sF_{r-}}
\bigl(\phi_r\bigl(X^{(n)}_{t,\cdot}\bigr)\bigr)\bigr)-G_r(X_{t,r-},\mE^{\sF_{r-}}(\phi
_r(X_{t,\cdot}))).
\]
As in (\ref{Eo1}), we can prove that for all $t\in[T,0]$,
%
\begin{equation}\label{Ep4}
\mE\Bigl(1_{A_n} \sup_{s\in[t,0]}\bigl|X^{(n)}_{t,s}-X_{t,s}\bigr|^\beta
\Bigr)\leq
C\mE\bigl(1_{A_n}\bigl|\xi^{(n)}-\xi\bigr|^\beta\bigr),
\end{equation}
where $C$ is independent of $n$.

Now, for any $\eps>0$, we have
\begin{eqnarray*}
P\Bigl\{\sup_{s\in[t,0]}\bigl|X^{(n)}_{t,s}-X_{t,s}\bigr|\geq\eps\Bigr\}
&\leq&
P\Bigl\{1_{A_n}\sup_{s\in[t,0]}\bigl|X^{(n)}_{t,s}-X_{t,s}\bigr|\geq\eps
\Bigr\}
+P(A^c_n)\\[-1pt]
&\leq&\frac{1}{\eps^\beta}\mE\Bigl(1_{A_n} \sup_{s\in
[t,0]}\bigl|X^{(n)}_{t,s}-X_{t,s}\bigr|^\beta\Bigr)+P(A^c_n)\\
&\leq&\frac{C}{\eps^\beta}\mE\bigl(1_{A_n}\bigl|\xi^{(n)}-\xi\bigr|^\beta\bigr)+P(A^c_n).
\end{eqnarray*}
The proof is then complete by letting $n\to\infty$.
\end{pf}
%
\begin{remark}
In the situation of Theorem~\ref{Th3}, the conclusion of this
proposition still holds, which can
be proven by the same procedure.
\end{remark}

The following lemma is a direct consequence of the uniqueness of solutions.

\begin{lemma}
Suppose that SFDE (\ref{SDE}) is uniquely solvable on the time interval
$(T,0]$. Then
for all $T<t_1<t_2<t_3\leq0$ and $\xi\in\sF_{t_1}$, we have
%
\begin{equation}\label{Ep1}
X_{t_2,t_3}(X_{t_1,t_2}(\xi))=X_{t_1,t_3}(\xi)\qquad \mbox{a.s.}
\end{equation}
Moreover, for any $T<t<s\leq0$, $x_i\in\mR^d, i=1,\ldots,n$ and
disjoint $\Lambda_i\in\sF_t, i=1,\ldots,n$
with $\bigcup_{i}\Lambda_i=\Omega$,
%
\begin{equation}\label{Ep2}
X_{t,s}\biggl(\sum_i1_{\Lambda_i} x_i\biggr)=\sum_i1_{\Lambda_i}X_{t,s}(x_i) \qquad\mbox{a.s.}
\end{equation}
\end{lemma}

\begin{pf}
For $T<t_1<t_2<s\leq0$, we can write
\[
X_{t_1,s}(\xi)=X_{t_1,t_2}(\xi)+\int_{(t_2,s]}G_r(X_{t_1,r-}(\xi
),\mE
^{\sF_{r-}}(\phi(X_{t_1,\cdot}(\xi))))\,\dif L_r \qquad\mbox{a.s.}
\]
On the other hand, if we set
\[
Y_s:=X_{t_2,s}(X_{t_1,t_2}(\xi))\qquad \forall s\in[t_2,0],
\]
then $Y_s$ satisfies
\[
Y_s=X_{t_1,t_2}(\xi)+\int_{(t_2,s]}G_r(Y_{r-},\mE^{\sF_{r-}}(\phi
(Y_\cdot)))\,\dif L_r\qquad \mbox{a.s.}
\]
Equality (\ref{Ep1}) follows by the uniqueness.

As for (\ref{Ep2}), noticing that for all $r\in(t,0]$,
\begin{eqnarray*}
&&\sum_i 1_{\Lambda_i}G_r(X_{t,r-}(x_i),\mE^{\sF_{r-}}(\phi
_r(X_{t,\cdot
}(x_i))))\\[-1pt]
&&\qquad =\sum_iG_r(1_{\Lambda_i}X_{t,r-}(x_i),1_{\Lambda_i}\mE
^{\sF_{r-}}
(\phi_r(X_{t,\cdot}(x_i))))\vadjust{\goodbreak}\\
& &\qquad=G_r\biggl(\sum_i 1_{\Lambda_i}X_{t,r-}(x_i),\mE^{\sF_{r-}}
\biggl(\phi_r\biggl(\sum_i 1_{\Lambda_i}X_{t,\cdot}(x_i)\biggr)\biggr)\biggr),
\end{eqnarray*}
it follows by the uniqueness as above.
\end{pf}

Now we can prove the following Markov property.
%
\begin{proposition}\label{Pro2}
In the situation of Theorem~\ref{Th1} or Theorem~\ref{Th3},
let $\{X_{t,s}(x); T\leq t<s\leq0\}$ be the solution family of SFDE
(\ref{SDE}) or (\ref{SDE2}).
Then for any $T\leq t_1<t_2<t_3\leq0$, $x\in\mR^d$ and bounded
continuous function $\varphi$,
we have
%
\begin{equation}\label{Es5}
\mE^{\sF_{t_2}}(\varphi(X_{t_1,t_3}(x)))=\mE(\varphi
(X_{t_2,t_3}(y)))|_{y=X_{t_1,t_2}(x)}\qquad   \mbox{a.s.}
\end{equation}
\end{proposition}

\begin{pf}
We only prove (\ref{Es5}) in the case of Theorem~\ref{Th1}.
By Proposition~\ref{Pro1}, the mapping $y\mapsto\mE(\varphi
(X_{t_2,t_3}(y))):=\Phi(y)$ is continuous.
So,\break $\Phi(X_{t_1,t_2}(x))$ is $\sF_{t_2}$-measurable. Thus, for proving
(\ref{Es5}),
it suffices to prove that for any $\Lambda\in\sF_{t_2}$,
\[
\mE(1_\Lambda\varphi(X_{t_1,t_3}(x)))=\mE(1_\Lambda\Phi(X_{t_1,t_2}(x))).
\]
Let $\xi^{(n)}=\sum_{i=1}^{m_n} x_i 1_{\Lambda_i}$ be a sequence of
simple functions, where $x_i\in\mR^d$,
$\Lambda_i\in\sF_{t_2}$ disjoint and $\bigcup_i\Lambda_i=\Omega$,
and such that
\[
\xi^{(n)}\to X_{t_1,t_2}(x) \qquad \mbox{in $L^\beta$ as }   n\to
\infty.
\]
By Proposition~\ref{Pro1} again, we have
\begin{eqnarray*}
\mE(1_\Lambda\varphi(X_{t_1,t_3}(x)))&\stackrel{\mathrm{\scriptsize(\ref{Ep1})}}{=}&
\mE(1_\Lambda\varphi(X_{t_2,t_3}(X_{t_1,t_2}(x))))\\
&=&\lim_{n\to\infty}\mE\bigl(1_\Lambda\varphi\bigl(X_{t_2,t_3}\bigl(\xi
^{(n)}\bigr)\bigr)\bigr)\\
&\stackrel{\scriptsize\mathrm{(\ref{Ep2})}}{=}&\lim_{n\to\infty}\sum_{i=1}^{m_n}\mE
(1_\Lambda1_{\Lambda_i}\varphi(X_{t_2,t_3}(x_i))).
\end{eqnarray*}
Since $X_{t_2,t_3}(x_i)$ is $\sF_{t_2,t_3}$-measurable and independent
of $\sF_{t_2}$, we further have
\begin{eqnarray*}
\mE(1_\Lambda\varphi(X_{t_1,t_3}(x)))
&=&\lim_{n\to\infty}\sum_{i=1}^{m_n}\mE(1_\Lambda1_{\Lambda
_i}\Phi
(x_i))\\
&=&\lim_{n\to\infty}\mE\bigl(1_\Lambda\Phi\bigl(\xi^{(n)}\bigr)\bigr)
=\mE(1_\Lambda\Phi(X_{t_1,t_2}(x))).
\end{eqnarray*}
The proof is complete.
\end{pf}

\subsection{Locally maximal solutions}\label{sec2.5}
Now, suppose that $\phi$ takes the following form:
%
\begin{equation}\label{Es6}
\phi_s(\omega)=\varphi(\omega(0))+\int^0_s f_r(\omega(r))\,\dif
r, \qquad \omega\in\mD,
\end{equation}
where $\varphi\dvtx \mR^d\to\mR^k$ and $f\dvtx \mR_-\times\mR^d\to\mR^k$ satisfy
that for some $K_2>0$
and all $s\in\mR_-$ and $x,x'\in\mR^d$,
%
\begin{equation}\label{EP2}
|\varphi(x)-\varphi(x')|+|f_s(x)-f_s(x')|\leq K_2|x-x'|.
\end{equation}
In this case, we have the following existence result of a unique
maximal solution.
%
\begin{theorem}\label{Th2}
Assume that (\ref{EP2}), \textup{(\textbf{H}$_G$)} and \textup{(\textbf{H}$^\beta_\nu$)} hold
for some $\beta>1$.
Then there exists a time $T=T(K_1,K_2,\sA,\beta)\in[-\infty,0)$ such
that SFDE (\ref{SDE})
is solvable on $(T,0]$ for any initial value $x\in\mR^d$, and if $T$ is
finite, then
%
\begin{equation}\label{Ep55}
\lim_{t\downarrow T}\|u_t\|_{\mathrm{Lip}}:=\sup_{x\not= x'\in\mR^d}\frac
{|u_t(x)-u_t(x')|}{|x-x'|}=+\infty,
\end{equation}
where
%
\begin{equation}\label{EP3}
u_t(x):=\mE\biggl(\varphi(X_{t,0}(x))+\int^0_t f_s(X_{t,s}(x))\,\dif
s\biggr).
\end{equation}
Moreover, the family of solutions $\{X_{t,s}(x),T<t<s\leq0,x\in\mR
^d\}
$ is unique in the class that
for all $T<t_1<t_2<t_3\leq0$ and $x\in\mR^d$,
\[
X_{t_1,t_2}(x)\in L^\beta(\Omega,\sF_{t_1,t_2},P),\qquad   X_{t_1,t_3}(x)=X_{t_2,t_3}(X_{t_1,t_2}(x))\qquad \mbox{a.s.}
\]
We also have the following uniform estimate: for any $T'\in(T,0)$ and
$x\in\mR^d$,
%
\begin{equation}\label{EL66}
\sup_{t\in[T',0]}\mE\Bigl(\sup_{s\in[t,0]}|X_{t,s}(x)|^\beta
\Bigr)\leq
C_{T',x}.
\end{equation}
\end{theorem}

\begin{pf}
First of all, let $T_1$ be the existence time in Theorem~\ref{Th1}.
By (\ref{Ep4}), there exists a constant $C=C(K_1,K_2,\sA,\beta)>0$
such that for all $x,x'\in\mR^d$ and $t\in[T_1,0]$,
\[
\mE\Bigl(\sup_{s\in[t,0]}|X_{t,s}(x)-X_{t,s}(x')|^\beta
\Bigr)\leq
C|x-x'|^\beta.
\]
Using this estimate and (\ref{EP2}), it is easy to check that
\[
\|u_{T_1}\|_{\mathrm{Lip}}<+\infty.
\]
Next, we consider the following SFDE on $[t,T_1]$:
\begin{eqnarray*}
X_{t,s}(x)&=&x+\int_{(t,s]}G_r\biggl(X_{t,r-}(x),\\
&&\hspace*{42pt}\qquad\mE^{\sF
_{r-}}\biggl(u_{T_1}(X_{t,T_1}(x))
+\int^{T_1}_r f_{r'}(X_{t,r'}(x))\,\dif r'\biggr)\biggr)\,\dif L_r.
\end{eqnarray*}
Repeating the proof of Theorem~\ref{Th1}, one can find another
$T_2<T_1$ so that this SFDE
is uniquely solvable on $[T_2,T_1]$. Meanwhile, one can patch up the
solution by setting
\[
X_{t,s}(x):=X_{T_1,s}(X_{t,T_1}(x)) \qquad \forall s\in[T_1,0],\ t\in[T_2,T_1].
\]
It is easy to verify that $\{X_{t,s}(x), T_2\leq t<s\leq0,x\in\mR^d\}$
solves SFDE (\ref{SDE}) on $[T_2,0]$.
Proceeding this construction, we obtain a sequence of times
\[
0>T_1>T_2>\cdots>T_n\downarrow T,
\]
and a family of solutions
\[
\{X_{t,s}(x), T<t<s\leq0,x\in\mR^d\}.
\]
From the construction of $T$, one knows that (\ref{Ep55}) holds. As for
the uniqueness, it can be proved
piecewisely on each $[T_n, T_{n-1}]$. Estimate (\ref{EL66}) follows
from (\ref{EL6}) and induction.
\end{pf}
%
\begin{remark}
By this theorem, for obtaining the global solution, it suffices to give
an a priori estimate
for $\|u_T\|_{\mathrm{Lip}}=\|\nabla u_T\|_\infty$.
\end{remark}

The following result can be proved similarly. We omit the details.
%
\begin{theorem}\label{Th22}
In addition to (\ref{EP2}) and \textup{(\textbf{H}$_G$)}, we assume that $G,
\varphi
$ and $f$
are uniformly bounded. Then there exists a time $T=T(K_1,K_2,\sA)\in
[-\infty,0)$ such that SFDE (\ref{SDE2})
is solvable on $(T,0]$ and estimate (\ref{Ep55}) holds provided
$T>-\infty$.
Moreover, the family of solutions $\{X_{t,s}(x),T<t<s\leq0,x\in\mR
^d\}
$ is unique in the class that
for all $T<t_1<t_2<t_3\leq0$ and $x\in\mR^d$,
\[
X_{t_1,t_2}(x)\in\sF_{t_1,t_2},\qquad  X_{t_1,t_3}(x)=X_{t_2,t_3}(X_{t_1,t_2}(x)) \qquad\mbox{a.s.}
\]
\end{theorem}

\section{Application to quasi-linear partial integro-differential equations}\label{sec3}

In this section, we establish the connection between stochastic
functional differential equations
and a class of quasi-linear partial integro-differential equations. For
this aim, we
consider $\phi$ taking the form of (\ref{Es6}) and assume that for some
$k\in\mN$,
(\textbf{H}$_k$), $G,f$ and $\varphi$ are continuous functions in $s,x,u$,
and for any $j=1,\ldots,k$, $\nabla^j G_s(x,u)$, $\nabla^j f_s(x)$,
$\nabla^j\varphi(x)$
are uniformly bounded continuous functions with respect to $s\in\mR_-$,
where $\nabla^j$ denotes the $j$th order
gradient with respect to $x,u$. We also denote
%
\begin{equation}\label{KK}
\sK:=\sup_{s\in\mR_-}(\|\nabla G_s\|_\infty+\|\nabla f_s\|
_\infty
)+\|\nabla\varphi\|_\infty.
\end{equation}

Under this assumption, it is clear that (\ref{EP2}) and (\textbf{H}$_G$) hold.
Let $u_t(x)$ be defined by (\ref{EP3}). By Theorem~\ref{Th2}, the mapping
$x\mapsto u_t(x)$ is Lipschitz continuous. However, it is in general
not $C^2$-differentiable
since we have poor integrabilities for $\nabla X_{t,s}(x)$. We shall
divide two cases to discuss this problem.

\subsection{\texorpdfstring{Unbounded data and $\nu$ has finite moments of arbitrary orders}
{Unbounded data and nu has finite moments of arbitrary orders}}\label{sec3.1}
In this subsection, we consider equation (\ref{SDE}), and assume that
(\textbf{H}$_k$) holds for some $k\geq3$,
and (\textbf{H}$^\beta_\nu$) holds for all $\beta\geq2$.
In this case, we can write
\[
L_t=\hat bt+W^A_t+\int_{\mR^m}z\tilde N(t,\dif z),
\]
where $\hat b:=b+\int_{|z|>1}z\nu(\dif z)\in\mR^m$.

Let $T<0$ be the maximal time given in Theorem~\ref{Th2} and $\{
X_{t,s}(x), T<t<s\leq0, x\in\mR^d\}$
the solution family of equation (\ref{SDE}). For simplicity of
notation, below we shall write
%
\begin{eqnarray}\label{EP0}
\cG_{t,r}&:=&\cG_{t,r}(x)
\nonumber
\\[-9pt]
\\[-9pt]
\nonumber
&:=&G_r\biggl(X_{t,r-}(x),\mE^{\sF
_{r-}}\biggl(\varphi(X_{t,0}(x))
+\int^0_r f_{r'}(X_{t,r'}(x))\,\dif r'\biggr)\biggr).
\end{eqnarray}
Let $g\dvtx \mR^d\to\mR^k$ be a $C^2$-function with bounded first and second
order partial derivatives.
By It\^o's formula (cf.~\cite{Ap}, page~226, Theorem~4.4.7), we have
%
\begin{eqnarray}\label{EL3}
g(X_{t,s})&=&g(x)+\int_{(t,s]} \int_{\mR^m}[g(X_{t,r-}+\cG
_{t,r}\cdot z)-g(X_{t,r-})
\nonumber\\[-2pt]
&&\hspace*{100pt}\qquad{}-\partial _i g(X_{t,r-})\cG^{ij}_{t,r} z_j]\nu(\dif z)\,\dif
r
\nonumber
\\[-9pt]
\\[-9pt]
\nonumber
&&{} +\int_{(t,s]}\partial _i g(X_{t,r-})\cG^{ij}_{t,r} \hat
b_j\,\dif r
\\[-2pt]
&&{}+\frac{1}{2}\int_{(t,s]}\partial _i\partial _jg(X_{t,r-})(\cG
^{\mathrm
{t}}_{t,r} A\cG
_{t,r})^{ij}\,\dif r+M^g_{t,s},\nonumber
\end{eqnarray}
where
\begin{eqnarray*}
M^g_{t,s}&:=&\int_{(t,s]} \int_{\mR^m}[g(X_{t,r-}+\cG
_{t,r}\cdot
z)-g(X_{t,r-})]\tilde N(\dif r,\dif z)
\\[-2pt]
&&{}+\int_{(t,s]}\partial _ig(X_{t,r-})\cG^{ij}_{t,r} \,\dif(W^A_r)^j
\end{eqnarray*}
is a square integrable ($\sF_{t,s}$)-martingale by (\ref{EL66}). Here
and below, the superscript ``t''
denotes the transpose of a matrix.

Fix $t\in(T,0]$ and $h>0$ so that $t-h\in(T,0]$. By taking expectations
for both sides of (\ref{EL3}), we have
\[
\frac{1}{h}[\mE g(X_{t-h,t})-g(x)]=I^g_1(h)+I^g_2(h)+I^g_3(h),
\]
where
\begin{eqnarray*}
I^g_1(h)&:=&\frac{1}{h}\mE\biggl(\int_{t-h}^t \int_{\mR^m}
[g(X_{t-h,r}+\cG_{t-h,r}\cdot z)-g(X_{t-h,r})
\\[-2pt]
&&\hspace*{92pt}\qquad{}-\partial _i g(X_{t-h,r})\cG^{ij}_{t-h,r} z_j]\nu(\dif z)\,\dif
r
\biggr),\\[-2pt]
I^g_2(h)&:=&\frac{1}{h}\mE\biggl(\int^t_{t-h}\partial
_ig(X_{t,r})\cG
^{ij}_{t-h,r} \hat b_j\,\dif r\biggr),\\[-2pt]
I^g_3(h)&:=&\frac{1}{2h}\mE\biggl(\int^t_{t-h}\partial _i\partial
_jg(X_{t,r})(\cG
^{\mathrm{t}}_{t-h,r} A\cG_{t-h,r})^{ij}\,\dif r\biggr).
\end{eqnarray*}
We have:
%
\begin{lemma}\label{Le2}
As $h\downarrow0$, it holds that
\begin{eqnarray*}
I^g_1(h)&\to&\int_{\mR^m}
\bigl[g\bigl(x+\sG_t(x)\cdot z\bigr)-g(x)-\partial _i g(x)\sG^{ij}_t(x)\cdot
z_j
\bigr]\nu
(\dif z),\\[-2pt]
I^g_2(h)&\to&\partial _ig(x)\sG^{ij}_t(x)\hat b_j,\qquad  I^g_3(h)\to
\tfrac
{1}{2}\partial
_i\partial _jg(x)(\sG^{\mathrm{t}}_t(x)A\sG_t(x))^{ij},
\end{eqnarray*}
where
\[
\sG_t(x):=G_t\biggl(x,\mE\biggl(\varphi(X_{t,0}(x))+\int^0_t
f_s(X_{t,s}(x))\,\dif s\biggr)\biggr).
\]
\end{lemma}

\begin{pf}
We only prove the first limit, the others are analogous. By the change
of variables, we can write
\begin{eqnarray*}
I^g_1(h)&=&\mE\biggl(\int^1_0   \int_{\mR^m}
[g(X_{t-h,t-hs}+\cG_{t-h,t-hs}\cdot z)\\[-2pt]
&&\hspace*{48pt}{}-g(X_{t-h,t-hs})
-\partial _i g(X_{t-h,t-hs})\cG^{ij}_{t-h,t-hs} z_j] \nu(\dif
z)\,\dif
s\biggr).
\end{eqnarray*}
Notice that
\begin{eqnarray*}
&&X_{t-h,t-hs}(x)-x\\[-2pt]
&&\qquad=\int_{(t-h,t-hs]} \int_{\mR^m}\cG
_{t-h,r}(x)\cdot
z\tilde N(\dif r,\dif z)
+\int_{(t-h,t-hs]}\cG_{t-h,r}(x)\cdot\hat b\,\dif r\\[-2pt]
&&\qquad\quad{} +\int_{(t-h,t-hs]}\cG_{t-h,r}(x)\,\dif W^A_r\\[-2pt]
&&\qquad=:J_1(h)+J_2(h)+J_3(h).
\end{eqnarray*}
By the isometric property of stochastic integrals, we have
\begin{eqnarray*}
\mE|J_1(h)|^2&=&\mE\biggl(\int_{t-h}^{t-hs} \int_{\mR^m}|\cG
_{t-h,r}(x)\cdot z|^2\nu(\dif z)\,\dif r\biggr)\\[-2pt]
&\leq&|h|\mE\biggl(\sup_{r\in[t-h,t-hs]}|\cG_{t-h,r}(x)|^2\int
_{\mR
^m}|z|^2\nu(\dif z)\biggr)\\[-2pt]
&\leq& C|h|\mE\Bigl(1+\sup_{r\in[t-h,0]}|X_{t-h,r}(x)|^2
\Bigr)\stackrel
{\scriptsize(\ref{EL66})}{\to}0\qquad\mbox{as $h\downarrow0$,}
\end{eqnarray*}
where in the last inequality, we used Doob's maximal inequality
and that $G,\phi$ and~$f$ in definition (\ref{EP0})
of $\cG_{t,r}$ are linear growth in $x$ and $u$, respectively.
Similarly,
\[
\mE|J_2(h)|^2+\mE|J_3(h)|^2\to0\qquad\mbox{as }h\downarrow0.
\]
Hence, for fixed $t,s,x$,
%
\begin{equation}\label{EL4}
\lim_{h\downarrow0}\mE|X_{t-h,t-hs}(x)-x|^2=0.
\end{equation}
Noticing that
%
\begin{equation}\label{For}
g(x+y)-g(x)=y\cdot\int^1_0\nabla g(x+\theta y)\,\dif\theta,
\end{equation}
we have
\begin{eqnarray*}
&&\mE|g(X_{t-h,t-hs}+\cG_{t-h,t-hs}\cdot z)-g(X_{t-h,t-hs})
-\partial _i g(X_{t-h,t-hs})\cG^{ij}_{t-h,t-hs} z_j|\\
&&\qquad=\mE\biggl|\biggl(\int^1_0[\partial _i g(X_{t-h,t-hs}+\theta\cG
_{t-h,t-hs}\cdot z)
-\partial _i g(X_{t-h,t-hs})]\,\dif\theta\biggr)\\
&&\hspace*{245pt}{}\times \cG^{ij}_{t-h,t-hs}
z_j\biggr|\\
&&\qquad\leq C\mE|\cG_{t-h,t-hs}|^2 |z|^2\leq
C\mE\Bigl(1+\sup_{r\in[t-h,0]}|X_{t-h,r}(x)|^2
\Bigr)|z|^2\stackrel
{\scriptsize(\ref{EL66})}{\leq} C|z|^2,
\end{eqnarray*}
where the second-to-last inequality is the same as above, and the
constant $C$ is independent of $h,s,z$.
Thus, for proving the first limit, by the dominated convergence
theorem, it suffices to prove that
for fixed $s\in[0,1]$ and $z\in\mR^m$,
\begin{eqnarray*}
&&\mE\biggl(\biggl(\int^1_0[\partial _i g(X_{t-h,t-hs}+\theta\cG
_{t-h,t-hs}\cdot z)
-\partial _i g(X_{t-h,t-hs})]\,\dif\theta\biggr)\cG^{ij}_{t-h,t-hs}
z_j\biggr)\\
&&\qquad \to\biggl(\int^1_0\bigl[\partial _i g\bigl(x+\theta\sG_t(x)\cdot z\bigr)
-\partial _i g(x)\bigr]\,\dif\theta\biggr)\sG^{ij}_t(x) z_j\qquad \mbox{as
$h\downarrow0$.}
\end{eqnarray*}
By (\ref{EL4}) and Remark~\ref{Re1}, this limit is easily obtained.
\end{pf}

We also need the following differentiability of the solution
$X_{t,s}(x)$ with respect to $x$ in the $L^p$-sense.
%
\begin{lemma}\label{Th4}
For any $p\geq2$, there exists a time $T_*=T_*(p,k,\sA,\sK)\in(T,0)$,
where $\sA$ is defined by (\ref{Ch}), and $\sK$ is defined by (\ref{KK}),
such that for any $T_*\leq t\leq s\leq0$,
the mapping $x\mapsto X_{t,s}(x)$ is $C^{k-1}$-differentiable in the
$L^p$-sense and for any $j=1,\ldots,k-1$,
\[
\sup_{x\in\mR^d}\sup_{s\in[t,0]}\mE|\nabla^j
X_{t,s}(x)|^p<+\infty.
\]
\end{lemma}

\begin{pf}
Since the proof is standard (cf.~\cite{Pr}, Theorem~39 or~\cite{Ku}, Section~4.6), we sketch it.
Let $\{e_i,i=1,\ldots, d\}$ be the canonical basis of $\mR^d$. For
$\delta>0$
and $i=1,\ldots, d$, define
\[
X^{\delta,i}_{t,s}:=X^{\delta,i}_{t,s}(x)=\frac{X_{t,s}(x+\delta
e_i)-X_{t,s}(x)}{\delta}
\]
and
\[
\cG^{\delta,i}_{t,s}:=\cG^{\delta,i}_{t,s}(x)=\frac{\cG
_{t,s}(x+\delta
e_i)-\cG_{t,s}(x)}{\delta},
\]
where $\cG_{t,s}(x)$ is defined by (\ref{EP0}). Then,
%
\begin{equation}\label{EU2}
\qquad X^{\delta,i}_{t,s}=e_i+\int_{(t,s]} \int_{\mR^m}\cG^{\delta
,i}_{t,r}\cdot z\tilde N(\dif r,\dif z)
+\int_{(t,s]}\cG^{\delta,i}_{t,r}\cdot\hat b\,\dif r+\int_{(t,s]}\cG
^{\delta,i}_{t,r}\,\dif W^A_r.
\end{equation}
As in (\ref{EP5}), by Burkholder's inequality, we have that for any
$p\geq2$,
%
\begin{equation}\label{EP1}
\mE\Bigl(\sup_{r\in[t,0]}|X^{\delta,i}_{t,r}|^p\Bigr)\leq
C_{p,\sA}|t|
\mE\Bigl(\sup_{r\in[t,0]}|\cG^{\delta,i}_{t,r}|^p\Bigr).
\end{equation}
Moreover, by (\textbf{H}$_k$) and Doob's maximal inequality, we easily
derive that
\[
\mE\Bigl(\sup_{r\in[t,0]}|\cG^{\delta,i}_{t,r}|^p\Bigr)
\leq C_{p,\sK}\mE\Bigl(\sup_{r\in[t,0]}|X^{\delta
,i}_{t,r}|^p\Bigr).
\]
Substituting this into (\ref{EP1}), we find that for some $C_{p,\sA
,\sK
}>0$ independent of $x,t$ and $\delta$,
\[
\mE\Bigl(\sup_{r\in[t,0]}|X^{\delta,i}_{t,r}|^p\Bigr)\leq
C_{p,\sA,\sK
}|t|\mE\Bigl(\sup_{r\in[t,0]}|X^{\delta,i}_{t,r}|^p\Bigr).
\]
From this, we deduce that there exists a time $T_*=T_*(p,\sA,\sK)\in
(T,0)$ such that for all $t\in[T_*,0]$,
%
\begin{equation}\label{EL2}
\sup_{\delta\in(0,1)}\sup_{x\in\mR^d}\sup_{t\in[T_*,0]}
\mE\Bigl(\sup_{r\in[t,0]}|X^{\delta,i}_{t,r}(x)|^p\Bigr)<+\infty
.
\end{equation}
On the other hand, let $Y^i_{t,s}=Y^i_{t,s}(x)$ satisfy the following SFDE:
%
\begin{eqnarray}\label{NN3}
Y^i_{t,s}&=&e_i+\int^s_t\nabla_x G_r\biggl(X_{t,r-}(x),\mE^{\sF
_{r-}}\biggl(\varphi(X_{t,0}(x))
+\int^0_r f_{r'}(X_{t,r'}(x))\,\dif r'\biggr)\biggr)\nonumber\\
&&\hspace*{35pt}{}\times Y^i_{t,r}\,\dif
L_r
\nonumber
\\[-8pt]
\\[-8pt]
\nonumber
&&{} +\int^s_t\nabla_y G_r(X_{t,r},\mE^{\sF_{r-}}(\phi
_r(X_{t,\cdot
})))\\
&&\hspace*{28pt}{}\times
\mE^{\sF_{r-}}\biggl(\nabla\varphi(X_{t,0})Y^i_{t,0}
+\int^0_r\nabla f_r'(X_{t,r'})Y^i_{t,r'}\,\dif r'\biggr)\,\dif
L_r,\nonumber
\end{eqnarray}
which can be solved on $[T_*,0]$ as in Theorem~\ref{Th1}.
Using the uniform estimate (\ref{EL2}) and formula (\ref{For}), it is
not hard to deduce that
\[
\lim_{\delta\to0}\mE\Bigl(\sup_{r\in[t,0]}|X^{\delta
,i}_{t,r}(x)-Y^i_{t,r}(x)|^p\Bigr)=0.
\]
In particular,
\[
\sup_{x\in\mR^d}\sup_{t\in[T_*,0]}\mE\Bigl(\sup_{r\in
[t,0]}|Y^i_{t,r}(x)|^p\Bigr)<+\infty.
\]
The higher derivatives can be estimated similarly from (\ref{NN3}).
\end{pf}

Now we can prove the following result, which was originally due to~\mbox{\cite{Bl,Ta,Fre}}.

\begin{theorem}\label{Th5}
Assume that \textup{(\textbf{H}$_k$)} holds for some $k\geq3$, and \textup{(\textbf{H}$^\beta_\nu$)} holds for all $\beta\geq2$.
Let $\{X_{t,s}(x), T<t\leq s\leq0,x\in\mR^d\}$ be the maximal solution
of SFDE~(\ref{SDE2}) in Theorem~\ref{Th2},
and $u_t(x)$ be defined by
%
\begin{equation}\label{EE}
u_t(x):=\mE\varphi(X_{t,0}(x))+\mE\biggl(\int^0_t
f_s(X_{t,s}(x))\,\dif
s\biggr).
\end{equation}
Then there exists a time $T_*=T_*(k,\sA,\sK)\in(T,0)$ such that for
each $t\in[T_*,0]$,
$x\mapsto u_t(x)$ has bounded derivatives up to ($k-1$)-order,
and solves the following quasi-linear partial integro-differential equation:
\[
u_t(x)=\varphi(x)+\int^0_t[\cL^{\mathrm{c}}u_s(x)+\cL^{\mathrm
{d}}u_s(x)+f_s(x)]\,\dif s\qquad
\forall(t,x)\in[T_*,0]\times\mR^d,
\]
where
\[
\cL^{\mathrm{c}}u_t(x):=\partial _iu_t(x)G^{ij}_t(x,u_t(x))\hat b_j
+\tfrac{1}{2}\partial _i\partial _ju_t(x)(G^{\mathrm
{t}}_t(x,u_t(x))AG_t(x,u_t(x)))^{ij}
\]
and
\[
\cL^{\mathrm{d}}u_t(x):=\int_{\mR^m}\bigl[u_t\bigl(x+G_t(x,u_t(x))\cdot z\bigr)-u_t(x)
-\partial _i u_t(x)G^{ij}_t(x,u_t(x))\cdot z_j\bigr]\nu(\dif z).
\]
\end{theorem}

\begin{pf}
We follow the argument of Friedman~\cite{Fr}. By Proposition \ref
{Pro2}, for $T<t-h<t\leq0$, we have
\begin{eqnarray*}
u_{t-h}(x)&=&\mE\bigl[(\mE\varphi
(X_{t,0}(y)))|_{y=X_{t-h,t}(x)}\bigr]
+\mE\biggl[\mE\biggl(\int^0_t f_r(X_{t,r}(y))\,\dif r\biggr)\bigg|_{y=
X_{t-h,t}(x)}\biggr]\\
&&{} +\mE\biggl(\int^t_{t-h} f_r(X_{t-h,r}(x))\,\dif r\biggr)\\
&=&\mE u_t(X_{t-h,t}(x))+\mE\biggl(\int^t_{t-h} f_r(X_{t-h,r}(x))\,\dif
r\biggr).
\end{eqnarray*}
By Lemma~\ref{Th4}, it is easy to see that there exists a time
$T_*=T_*(k,\sA,\sK)<0$ such that for each $t\in[T_*,0]$,
$u_t(x)$ has bounded derivatives\vadjust{\goodbreak} up to $(k-1)$-order.
Thus, we can invoke Lemma~\ref{Le2} to derive that
\begin{eqnarray*}
&&\frac{1}{h}\bigl(u_{t-h}(x)-u_t(x)\bigr)\\
&&\qquad=\frac{1}{h}\bigl(\mE u_t(X_{t-h,t}(x))-u_t(x)\bigr)
+\frac{1}{h}\mE\biggl(\int^t_{t-h} f_r(X_{t-h,r}(x))\,\dif r\biggr)\\
&&\qquad\to\cL^{\mathrm{c}}u_t(x)+\cL^{\mathrm{d}}u_t(x)+f_t(x)\qquad\mbox{as
$h\downarrow0$}.
\end{eqnarray*}
On the other hand, from the above proof, it is also easy to see that
for fixed $x\in\mR^d$,
$t\mapsto u_t(x)$ is Lipschitz continuous. Hence,
\[
u_t(x)-\varphi(x)=-\int^0_t\partial _s u_s(x)\,\dif s=\int^0_t[\cL
^{\mathrm
{c}}u_s(x)+\cL^{\mathrm{d}}u_s(x)+f_s(x)]\,\dif s.
\]
The proof is thus complete.
\end{pf}

\subsection{Bounded data and constant big jump}\label{sec3.2}

In this subsection we assume that (\textbf{H}$_k$) holds for some $k\geq
3$, and $G,\varphi$ and $f$ are uniformly bounded
and continuous functions. Consider the following SFDE:
\begin{eqnarray*}
X_{t,s}(x)&=&x+\int_{(t,s]} \int_{|z|<1}\cG_{t,r}(x)\cdot z\tilde
N(\dif r,\dif z)
+\int_{(t,s]} \int_{|z|\geq1}z N(\dif r,\dif z)\\
&&{} +\int_{(t,s]}\cG_{t,r}(x)\cdot b\,\dif r+\int_{(t,s]}\cG
_{t,r}(x)\,\dif W^A_r,
\end{eqnarray*}
where $\cG_{t,r}(x)$ is defined by (\ref{EP0}).
In this case, Lemmas~\ref{Le2} and~\ref{Th4} still hold. We just want
to mention that (\ref{EL4})
should be replaced by
\[
X_{t-h,t-hs}(x)\to x \qquad\mbox{in probability as $h\downarrow0$,}
\]
and (\ref{EU2}) becomes
\[
X^{\delta,i}_{t,s}=e_i+\int_{(t,s]} \int_{|z|\leq1}\cG^{\delta
,i}_{t,r}\cdot z\tilde N(\dif r,\dif z)
+\int_{(t,s]}\cG^{\delta,i}_{t,r}\cdot b\,\dif r+\int_{(t,s]}\cG
^{\delta
,i}_{t,r}\,\dif W^A_r.
\]

Thus, the following result can be proved along the same lines as in
Theorem~\ref{Th5}. We omit the details.
%
\begin{theorem}\label{Th8}
Assume that \textup{(\textbf{H}$_k$)} holds for some $k\geq3$, and $G,\varphi$ and
$f$ are uniformly bounded
and continuous functions. Let $\{X_{t,s}(x), T<t\leq s\leq0,x\in\mR
^d\}$
be the short time solution of SFDE (\ref{SDE2}) in Theorem~\ref{Th22},
and $u_t(x)$ be defined by (\ref{EE}).
Then there exists a time $T_*=T_*(k,\sA,\sK)\in(T,0)$ such that for
each $t\in[T_*,0]$,
$x\mapsto u_t(x)$ has bounded derivatives up to ($k-1$)-order, and
solves the following quasi-linear partial integro-differential equation:
\[
u_t(x)=\varphi(x)+\int^0_t[\cL^{\mathrm{c}}u_s(x)+\cL^{\mathrm
{d}}u_s(x)+f_s(x)]\,\dif s\qquad
\forall(t,x)\in[T_*,0]\times\mR^d,
\]
where
\[
\cL^{\mathrm{c}}u_t(x):=\partial _iu_t(x)G^{ij}_t(x,u_t(x))b_j
+\tfrac{1}{2}\partial _i\partial _ju_t(x)(G^{\mathrm
{t}}_t(x,u_t(x))AG_t(x,u_t(x)))^{ij}
\]
and
\begin{eqnarray*}
\cL^{\mathrm{d}}u_t(x)&:=&\int_{|z|<1}\bigl[u_t\bigl(x+G_t(x,u_t(x))\cdot z\bigr)-u_t(x)
-\partial _i u_t(x)G^{ij}_t(x,u_t(x))\cdot z_j\bigr]\\
&&\hspace*{23pt}{}\times \nu(\dif z)\\
&&{} +\int_{|z|\geq1}[u_t(x+z)-u_t(x)]\nu(\dif z).
\end{eqnarray*}
\end{theorem}

\section{Semi-linear partial integro-differential equation: Existence
and uniqueness of weak solutions}\label{sec4}

In this section we consider the following semi-linear partial
integro-differential equation:
%
\begin{equation}\label{EE1}
\partial _t u_t+\cL_0 u_t+G^i_t(x, u_t)\partial _i
u_t+F_t(x,u_t)=0,\qquad  u_0=\varphi, t\leq0,
\end{equation}
where $\cL_0$ is the generator of the L\'evy process $L_t$ given by
(\ref{Ge}), and
%
\begin{eqnarray}\label{EE3}
\qquad G&\in&\cB\bigl(\mR_-; \mW^{1,\infty}(\mR^d\times\mR^k;\mR^d)\bigr),\qquad
F\in\cB\bigl(\mR_-; \mW^{1,\infty}(\mR^d\times\mR^k;\mR^k)\bigr),
\nonumber
\\[-8pt]
\\[-8pt]
\nonumber
\varphi&\in&\mW^{1,\infty}(\mR^d;\mR^k).
\end{eqnarray}
Here and below, $\mW^{1,\infty}$ denotes the space of bounded and
Lipschitz continuous functions,
$\cB$ or $\cB_{\mathrm{loc}}$ denotes the space of uniformly or locally bounded
measurable functions.

Let us first give the following definition about the maximal weak
solution for equation (\ref{EE1}).
%
\begin{definition}\label{EE4}
For $T\in[-\infty,0)$, we call $u\in\cB_{\mathrm{loc}}((T,0]; \mW^{1,\infty
}(\mR
^d;\mR^k))$ a maximal weak solution of
equation (\ref{EE1}) if
%
\begin{equation}\label{Max}
\lim_{t\downarrow T}\|\nabla u_t(x)\|_\infty=+\infty\qquad\mbox{when }  T>-\infty,
\end{equation}
and for all $\psi\in C^\infty_0(\mR^d;\mR^k)$ and $t\in(T,0]$,
%
\begin{equation}\label{EP9}
\qquad\langle u_t,\psi\rangle=\langle\varphi,\psi\rangle+\int
^0_t\langle
u_s,\cL_0^*\psi\rangle\,\dif r
+\int^0_t\langle G^i_s(u_s)\partial _iu_s+F_s(u_s),\psi\rangle\,\dif
s,
\end{equation}
where $\langle\varphi,\psi\rangle:=\int_{\mR^d}\varphi(x)\cdot
\psi
(x)\,\dif x$, and $\cL^*_0$ is the adjoint
operator of $\cL_0$ and given by
\[
\cL^*_0 \psi(x):=\frac{1}{2}a_{ij}\partial _i\partial _j \psi
-b_i\partial _i \psi
+\int_{\mR
^d}\bigl[\psi(x-z)-\psi(x)+1_{|z|<1}
\partial _i \psi(x)z_i\bigr]\nu(\dif z).
\]
\end{definition}

The main aim of this section is to prove the following existence and
uniqueness of a maximal weak solution
as well as the global solution for equation (\ref{EE1}).\vadjust{\goodbreak}
%
\begin{theorem}\label{Main}
\textup{(i)} (Locally maximal weak solution) Under (\ref{EE3}), there exists a
unique maximal weak solution $u_t(x)$ for equation (\ref{EE1})
in the sense of Definition~\ref{EE4}. Moreover, let $T$ be the maximal
existence time, then for any $t\in(T,0]$,
%
\begin{equation}\label{AS}
\|u_t\|_\infty\leq\|\varphi\|_\infty+|t|\sup_{s\in[t,0]}\|F_s\|
_\infty
.
\end{equation}\vspace*{-6pt}
\begin{longlist}[(iii)]
\item[(ii)] (Nonnegative solution) If for some $j=1,\ldots,k$, the components
$\varphi^j$ and $F^j$ are nonnegative, then
the corresponding component $u^j$ of weak solution in \textup{(i)} are also nonnegative.

\item[(iii)] (Global solution) Let $\Psi(\xi)$ be the L\'evy symbol defined in
(\ref{Lp9})
with $b=A=0$. If for some $\alpha\in(1,2)$,
%
\begin{equation}\label{EU3}
\mathrm{Re}(\Psi(\xi))\asymp|\xi|^\alpha \qquad\mbox{as } |\xi|\to
\infty
,
\end{equation}
where $a\asymp b$ means that for some $c_1,c_2>0$, $c_1 b\leq a\leq c_2
b$, then the maximal
existence time $T$ in \textup{(i)} equals to $-\infty$. In the case that $b=\nu
=0$ and $A$ is strictly positive,
then $T$ also equals to $-\infty$.
\end{longlist}
\end{theorem}

\begin{remark}
Since we have estimate (\ref{AS}), it is easy to see that the
assumption on $G$ in
(\ref{EE3}) can be replaced by
\[
G\in\cB\bigl(\mR_-; \mW^{1,\infty}(\mR^d\times\mB_R;\mR^d)\bigr)\qquad \forall R>0,
\]
where $\mB_R:=\{x\in\mR^k\dvtx  |x|\leq R\}$.
\end{remark}

For proving this theorem, let us begin with studying:

\subsection{Linear partial integro-differential equation}\label{sec4.1}
In this subsection, we firstly study the existence and uniqueness of
weak solutions for the following linear PIDE:
%
\begin{equation}\label{EE11}\quad
\partial _t u_t+\cL_0 u_t+G^i_t(x)\partial _i
u_t+H_t(x)u_t+f_t(x)=0,\qquad  u_0=\varphi
,\ t\leq0,
\end{equation}
where $G\dvtx \mR_-\times\mR^d\to\mR^d$, $H\dvtx \mR_-\times\mR^d\to\mR
^k\times\mR^k$,
$f\dvtx \mR_-\times\mR^d\to\mR^k$ and $\varphi\dvtx \mR^d\to\mR^k$ are bounded
measurable functions.

Let us start with the following case of smooth coefficients, which is
the~classical Feynman--Kac formula.
Here, the main point is to prove the uniqueness.

\begin{theorem}[(Feynman--Kac formula)]\label{EE7}
Assume that
\begin{eqnarray*}
G&\in&\cB(\mR_-; C^\infty_b(\mR^d;\mR^d)),\qquad
H\in\cB\bigl(\mR_-; C^\infty_b(\mR^d;\mR^k\times\mR^k)\bigr),
\\
f&\in&\cB(\mR_-; C^\infty_b(\mR^d;\mR^k)),\qquad
\varphi\in C^\infty_b(\mR^d;\mR^k),
\end{eqnarray*}
where $C^\infty_b$ denotes the space of bounded smooth functions with
bounded derivatives of all orders.
Let $\{X_{t,s}(x),t\leq s\leq0,x\in\mR^d\}$ solve the following SDE:
\[
X_{t,s}(x)=x+\int^s_t G_r(X_{t,r}(x))\,\dif r+\int^s_t\,\dif L_r,\vadjust{\goodbreak}
\]
and $\{Z_{t,s}(x),t\leq s\leq0,x\in\mR^d\}$ solve the following ODE:
\[
Z_{t,s}(x)=\mI_{m\times m}+\int^s_t H_r(X_{t,r}(x))\cdot
Z_{t,r}(x)\,\dif x.
\]
Define
%
\begin{equation}\label{EE5}
u_t(x):=\mE[Z_{t,0}(x)\varphi(X_{t,0}(x))]+\mE
\biggl[\int^0_t
Z_{t,r}(x)f_r(X_{t,r}(x))\,\dif r\biggr].
\end{equation}
Then $u\in C(\mR_-; C^\infty_b(\mR^d;\mR^k))$ uniquely solves the
following linear PIDE:
%
\begin{eqnarray}\label{EE111}
u_t(x)=\varphi(x)+\int^0_t[\cL_0 u_s(x)+G^i_s(x)\partial _i
u_s(x)+H_s(x)u_s(x)+f_s(x)]\,\dif s
\nonumber
\\[-8pt]
\\[-8pt]
\eqntext{ \forall(t,x)\in\mR_-\times\mR^d.}
\end{eqnarray}
\end{theorem}

\begin{pf}
By smoothing the time variable and then taking limits, as in Section~\ref{sec3},
by careful calculations,
one can find that $u$ defined by (\ref{EE5}) belongs to $C(\mR_-;
C^\infty_b(\mR^d;\mR^k))$ and satisfies (\ref{EE111}); see~\cite{Fr}, page~148, Theorem 5.3 \mbox{or~\cite{Bu-Fl-Ro,Zh2}}.

We now prove the uniqueness by the duality argument.
Let $\hat X_{t,s}(x)$ solve the following SDE:
\[
\hat X_{t,s}(x)=x-\int^s_t G_r(\hat X_{t,r}(x))\,\dif r-\int^s_t\,\dif L_r,
\]
and $\hat Z_{t,s}(x)$ solve the following ODE:
\[
\hat Z_{t,s}(x)=\mI_{m\times m}+\int^s_t [H_r(\hat
X_{t,r}(x))^{\mathrm{t}}
+\div G_r(\hat X_{t,r}(x))\mI_{m\times m}]\cdot\hat
Z_{t,r}(x)\,\dif x.
\]
Fix $T<0$ and $\psi\in C^\infty_0(\mR^d;\mR^k)$ define
\[
\hat u_t(x):=\mE[\hat Z_{T-t,0}(x)\psi(\hat X_{T-t,0}(x))].
\]
As above, one can check that
\[
\hat u_t(x)=\psi(x)+\int^t_T[\cL^*_0 \hat u_s(x)-G^i_s(x)\partial
_i\hat u_s(x)
+H_s(x)^{\mathrm{t}}\hat u_s(x)+\div G_s(x)\hat u_s(x)]\,\dif s.
\]
Let $u\in C(\mR_-; C^\infty_b(\mR^d;\mR^k))$ satisfy equation (\ref
{EE111}) with $\varphi=f=0$.
Then by the integration by parts formula, we have for almost all $t\in[T,0]$,
\begin{eqnarray*}
\partial _t\langle u_t,\hat u_t\rangle&=&-\langle\cL_0
u_t+G^i_t\partial _i
u_t+H_tu_t,\hat u_t\rangle
\\
&&{}+\langle u_t,\cL_0^*\hat u_t-G_t\partial _i\hat u_t+H^{\mathrm
{t}}_t\hat
u+\div G_t\hat u_t\rangle\\
&=&0.
\end{eqnarray*}
From this, we get
\[
\langle u_T,\psi\rangle=\langle u_T,\hat u_T\rangle=\langle u_0,\hat
u_0\rangle=0,
\]
which leads to $u_T(x)=0$ by the arbitrariness of $\psi$.\vadjust{\goodbreak}
\end{pf}

For $\ell\in\mN$, we introduce a family of mollifiers in $\mR^\ell$.
Let $\rho\dvtx \mR^\ell\to[0,1]$ be a smooth function satisfying that
\[
\rho(x)=0\qquad \forall|x|>1,\qquad  \int_{\mR^\ell}\rho(x)\,\dif x=1.
\]
We shall call $\{\rho_\eps(x):=\eps^{-\ell}\rho(x/\eps),\eps\in
(0,1)\}$
a family of mollifiers in $\mR^\ell$.

Next, we relax the regularity assumptions on $G,H,f$ and $\varphi$, and
prove the following:
%
\begin{theorem}\label{Th6}
Assume that
\begin{eqnarray*}
 G&\in&\cB(\mR_-; \mW^{1,\infty}(\mR^d;\mR^d)),\qquad
H\in\cB(\mR_-\times\mR^d;\mR^k\times\mR^k),
\\
f&\in&\cB(\mR_-; \mW^{1,\infty}(\mR^d;\mR^k)),\qquad
\varphi\in\mW^{1,\infty}(\mR^d;\mR^k).
\end{eqnarray*}
Let $u_t(x)$ be defined as in (\ref{EE5}). Then $u_t(x)\in\cB
_{\mathrm{loc}}(\mR
_-; \mW^{1,\infty}(\mR^d;\mR^k))$
is a unique weak solution of equation (\ref{EE11}) in the sense of
Definition~\ref{EE4}.
\end{theorem}

\begin{pf}
We only prove the uniqueness. As for the existence, it follows by
smoothing the coefficients and then taking limits
as done in Theorem~\ref{Th7} below.

Suppose that $u\in\cB_{\mathrm{loc}}(\mR_-; \mW^{1,\infty}(\mR^d;\mR^k))$
is a
weak solution of equation (\ref{EE11})
with $\varphi=f=0$ in the sense of Definition~\ref{EE4}. We want to
prove that $u\equiv0$.
Let $\rho_\eps$ be a family of mollifiers in $\mR^d$. Define
\[
u^\eps_t(x):=u_t*\rho_\eps(x),\qquad  G^\eps_t(x):=G_t*\rho_\eps(x),\qquad  H^\eps_t(x):=H_t*\rho_\eps(x).
\]
Taking $\psi(\cdot)=\rho_\eps(x-\cdot)$ in (\ref{EP9}), one finds that
$u^\eps_t(x)$ satisfies
\[
u^\eps_t(x)=\int^0_t[\cL_0 u^\eps_s(x)+G^{\eps,i}_s(x)\partial _i
u^\eps_s(x)
+H^\eps_s(x)u^\eps_s(x)+f^\eps_s(x)]\,\dif s,
\]
where
\[
f^\eps_s(x)=(G^i_s\partial _i u)*\rho_\eps(x)-G^{\eps
,i}_s(x)\partial _i u^\eps_s(x)
+(H_su_s)*\rho_\eps(x)-H^\eps_s(x)u^\eps_s(x).
\]
By the property of convolutions, we have
%
\begin{equation}
\|f^\eps_s\|_\infty\leq2\|G_s\|_\infty\|\nabla u\|_\infty+2\|H_s\|
_\infty\|u\|_\infty\label{UU0},
\end{equation}
and for fixed $s$ and Lebesgue almost all $x\in\mR^d$,
%
\begin{equation}\label{UU}
f^\eps_s(x)\to0, \qquad \eps\to0.
\end{equation}
Let $X^\eps_{t,s}(x)$ solve the following SDE:
%
\begin{equation}\label{NN2}
X^\eps_{t,s}(x)=x+\int^s_t G^\eps_r(X^\eps_{t,r}(x))\,\dif r+\int
^s_t\,\dif
L_r,
\end{equation}
and $Z^\eps_{t,s}(x)$ solve the following ODE:
%
\begin{equation}\label{NN1}
Z^\eps_{t,s}(x)=\mI_{m\times m}+\int^s_t H^\eps_r(X^\eps
_{t,r}(x))\cdot
Z^\eps_{t,r}(x)\,\dif x.
\end{equation}
By Theorem~\ref{EE7}, $u^\eps_t(x)$ can be uniquely represented by
\[
u^\eps_t(x):=\mE\biggl[\int^0_t Z^\eps_{t,s}(x)f^\eps_s(X^\eps
_{t,s}(x))\,\dif s\biggr].
\]
For completing the proof, it suffices to prove that for each $(t,x)\in
\mR_-\times\mR^d$,
\[
u^\eps_t(x)\to0, \qquad  \eps\to0.
\]
Since by (\ref{NN1}), $Z^\eps_{t,s}(x)$ is uniformly bounded with
respect to $x,\eps$ and $s\in[t,0]$,
we need only to show that for any nonnegative $\psi\in C^\infty
_0(\mR^d)$,
\[
I^\eps_t:=\mE\biggl[\int^0_t    \int_{\mR^d}|f^\eps_s|(X^\eps
_{t,s}(x))\psi(x)\,\dif x\,\dif s\biggr]\to0
,  \qquad \eps\to0.
\]
For any $R>0$, by (\ref{UU0}) and the change of variables, we have
%
\begin{eqnarray}\label{UU1}
I^\eps_t&\leq&\mE\biggl[\int^0_t    \int_{|X^\eps_{t,s}(x)|\leq R}
|f^\eps_s|(X^\eps_{t,s}(x))\psi(x)\,\dif x\,\dif s\biggr]\nonumber\\
&&{}+C_t\mE\biggl[\int^0_t    \int_{|X^\eps_{t,s}(x)|>R}\psi(x)\,\dif
x\,\dif
s\biggr]
\nonumber
\\[-8pt]
\\[-8pt]
\nonumber
&\leq&\mE\biggl[\int^0_t    \int_{B_R}
|f^\eps_s|(x)\psi(X^{\eps,-1}_{t,s}(x))\det(\nabla X^{\eps
,-1}_{t,s}(x))\,\dif x\,\dif s\biggr] \\
&&{} +C_t\int^0_t    \int_{\mR^d}P\{|X^\eps
_{t,s}(x)|>R\}
\psi(x)\,\dif x\,\dif s,\nonumber
\end{eqnarray}
where $X^{\eps,-1}_{t,s}(x)$ denotes the inverse of $x\mapsto X^{\eps
}_{t,s}(x)$.
From equation (\ref{NN2}), it is by now standard to prove that (e.g.,
see Kunita~\cite{Ku}, Lemma~4.3.1)
\[
\det(\nabla X^{\eps,-1}_{t,s}(x))=\exp\biggl\{-\int^s_t(\div G^\eps
_r)(X^\eps_{t,r}(X^{\eps,-1}_{t,s}(x)))\,\dif r\biggr\},
\]
which then yields
\[
C_0:=\sup_{\eps\in(0,1)}\sup_{x\in\mR^d}\sup_{s\in[t,0]}|\det
(\nabla
X^{\eps,-1}_{t,s}(x))|<+\infty.
\]
Thus, for fixed $R>0$, the first term in (\ref{UU1}) is less than
\[
C_0\|\psi\|_\infty\int^0_t    \int_{B_R}|f^\eps_s|(x)\,\dif x\,\dif
s\stackrel{\scriptsize(\ref{UU})}{\to} 0,  \qquad\eps\to0.
\]
Moreover, by equation (\ref{NN2}), we also have
\[
\lim_{R\to\infty}\sup_{\eps}P\{|X^\eps_{t,s}(x)|>R\}
\leq\lim_{R\to\infty}P\biggl\{|x|+\int^0_t\|G_r\|_\infty\,\dif
r+|L_s-L_t|>R\biggr\}=0.
\]
The proof is complete by first letting $\eps\to0$ and then $R\to
\infty
$ in (\ref{UU1}).
\end{pf}

As an easy corollary of this theorem, we first establish the uniqueness
for equation (\ref{EE1}).
%
\begin{theorem}
Under (\ref{EE3}), there exists, at most, one weak solution for
equation (\ref{EE1}).
\end{theorem}

\begin{pf}
Let $u^{(i)}\in\cB_{\mathrm{loc}}((T,0]; \mW^{1,\infty}(\mR^d;\mR^k)), i=1,2$
be two weak solutions of equation (\ref{EE1})
in the sense of Definition~\ref{EE4}. Define
$u_t(x):=u^{(1)}_t(x)-u^{(2)}_t(x)$.
Then $u_t(x)$ satisfies that for all $\psi\in C^\infty_0(\mR^d;\mR^k)$,
\[
\langle u_t,\psi\rangle=\int^0_t\langle u_s,\cL_0^*\psi\rangle
\,\dif
r+\int^0_t\bigl\langle G^i_s\bigl(u^{(1)}_s\bigr)\partial _iu_s,\psi\bigr\rangle\,\dif
s+\int
^0_t\langle H_su_s,\psi\rangle\,\dif s,
\]
where
\begin{eqnarray*}
H_s(x)&:=&\biggl(\int^1_0\nabla_uG^i_s\bigl(x,u^{(1)}_s(x)+\theta
u_s(x)\bigr)\,\dif
\theta\biggr)\partial _i u^{(2)}_s(x)
\\
&&{}+\int^1_0\nabla_uF_s\bigl(x,u^{(1)}_s(x)+\theta u_s(x)\bigr)\,\dif\theta.
\end{eqnarray*}
By (\ref{EE3}), it is easy to verify that
\[
(s,x)\mapsto G_s\bigl(x,u^{(1)}_s(x)\bigr)\in\cB_{\mathrm{loc}}((T,0];\mW^{1,\infty
}(\mR
^d;\mR^d))
\]
and
\[
(s,x)\mapsto H_s(x)\in\cB_{\mathrm{loc}}((T,0];\cB(\mR^d;\mR^k)).
\]
Thus, by Theorem~\ref{Th6}, we conclude that $u_t(x)=0$.
\end{pf}

\subsection{A special form: $F_t(x,u)=f_t(x)$ independent of $u$}\label{sec4.2}

Consider the following SFDE:
%
\begin{eqnarray}\label{SFEE}
 X_{t,s}(x)&=&x+\int^s_tG_r\biggl(X_{t,r}(x),\mE^{\sF_r}\biggl(\varphi(X_{t,0}(x))
-\int^0_r f_{r'}(X_{t,r'}(x))\,\dif r'\biggr)\biggr)\,\dif r
\nonumber\hspace*{-35pt}
\\[-8pt]
\\[-8pt]
\nonumber
&&{}+\int^s_t\,\dif L_r,\hspace*{-35pt}
\end{eqnarray}
where $G\dvtx \mR_-\times\mR^d\times\mR^k\to\mR^d$, $f\dvtx \mR_-\times
\mR^d\to\mR
^k$ and $\varphi\dvtx \mR^d\to\mR^k$
are bounded measurable functions.

We need the following continuous dependence of the solutions with
respect to the coefficients.
%
\begin{proposition}\label{Pro3}
Suppose that $(G^{(i)}, f^{(i)},\varphi^{(i)}), i=1,2$ are two groups
of bounded measurable functions, and
for some $K>0$ and all $t\in\mR_-$, $x,x'\in\mR^d$ and $u,u'\in\mR^k$,
\begin{eqnarray*}
&&\bigl|G^{(i)}_t(x,u)-G^{(i)}_t(x',u')\bigr|+\bigl|f^{(i)}_t(x)-f^{(i)}_t(x')\bigr|+\bigl|\varphi
^{(i)}(x)-\varphi^{(i)}(x')\bigr|\\
&&\qquad\leq K(|x-x'|+|u-u'|).
\end{eqnarray*}
Then there exists a time $T<0$ depending only on $K$ such that for all
$t\in[T,0]$ and $x,y\in\mR^d$,
%
\begin{eqnarray}
&&\sup_{s\in[t,0]}\mE\bigl|X^{(1)}_{t,s}(x)-X^{(2)}_{t,s}(y)\bigr|\nonumber\\
&&\qquad\leq2|x-y|
+2\int^0_t\bigl\|G^{(1)}_r-G^{(2)}_r\bigr\|_\infty\,\dif r
\\
&&\qquad\quad{} +2K|T|\biggl(\bigl\|\varphi^{(1)}-\varphi^{(2)}\bigr\|_\infty+\int^0_t
\bigl\|
f^{(1)}_r-f^{(2)}_r\bigr\|_\infty\,\dif r\biggr),\nonumber
\end{eqnarray}
where $X^{(i)}_{t,s}(x)$ is the solution of (\ref{SFEE}) corresponding
to $(G^{(i)}, f^{(i)},\varphi^{(i)})$.
\end{proposition}

\begin{pf}
Set
\[
Z_{t,s}:=X^{(1)}_{t,s}(x)-X^{(2)}_{t,s}(y).
\]
By (\ref{SFEE}) and the assumptions, we have
\begin{eqnarray*}
\mE|Z_{t,s}|&\leq&|x-y|+\int^s_t\bigl\|G^{(1)}_r-G^{(2)}_r\bigr\|_\infty\,\dif r
\\
&&{}+K\int^s_t\biggl(\mE|Z_{t,r}|+\bigl\|\varphi^{(1)}-\varphi^{(2)}\bigr\|
_\infty+K\mE
|Z_{t,0}|\\
&&\hspace*{22pt}\qquad{}+\int^0_r \bigl\|f^{(1)}_{r'}-f^{(2)}_{r'}\bigr\|_\infty\,\dif r'
+K\int^0_r\mE|Z_{t,r'}|\,\dif r'\biggr)\,\dif r\\
&\leq&|x-y|+\int^0_t\bigl\|G^{(1)}_r-G^{(2)}_r\bigr\|_\infty\,\dif r\\
&&{}+K|t|
\biggl(\bigl\|\varphi^{(1)}-\varphi^{(2)}\bigr\|_\infty+\int^0_t \bigl\|
f^{(1)}_r-f^{(2)}_r\bigr\|_\infty\,\dif r\biggr)\\
&&{} +K\int^s_t\mE|Z_{t,r}|\,\dif r+K^2|t|\biggl(\mE|Z_{t,0}|
+\int^0_t\mE|Z_{t,r}|\,\dif r\biggr)\\
&\leq&|x-y|+\int^0_t\bigl\|G^{(1)}_r-G^{(2)}_r\bigr\|_\infty\,\dif r\\
&&{}+K|t|
\biggl(\bigl\|\varphi^{(1)}-\varphi^{(2)}\bigr\|_\infty+\int^0_t \bigl\|
f^{(1)}_r-f^{(2)}_r\bigr\|_\infty\,\dif r\biggr)\\
&&{} +(K|t|+K^2|t|+K^2|t|^2)\sup_{r\in[t,0]}\mE|Z_{t,r}|.
\end{eqnarray*}
From this, we immediately conclude the proof.
\end{pf}
%
\begin{theorem}\label{Th7}
Assume that $(G,f,\varphi)$ are bounded measurable functions and satisfy
for some $K>0$ and all $t\in\mR_-$, $x,x'\in\mR^d$ and $u,u'\in\mR^k$,
\[
|G_t(x,u)-G_t(x',u')|+|f_t(x)-f_t(x')|+|\varphi(x)-\varphi(x')|
\leq K(|x-x'|+|u-u'|).
\]
Then there exists a time $T=T(K)<0$ such that
%
\begin{equation}\label{EU}
u_t(x):=\mE\varphi(X_{t,0}(x))+\mE\biggl(\int^0_t
f_r(X_{t,r}(x))\,\dif
r\biggr)
\end{equation}
is a unique weak solution of equation (\ref{EE1}) on $[T,0]$ in the
sense of Definition~\ref{EE4}.
\end{theorem}

\begin{pf}
Let ($G^\eps,f^\eps,\varphi^\eps$) be the smooth approximation of
($G,f,\varphi$) defined by
\begin{eqnarray*}
G^\eps_t(x,u)&:=&G*\rho^{(1)}_\eps(t,x,u), \qquad  f^\eps_t(x):=f*\rho
^{(1)}_\eps(t,x),\\
  \varphi^\eps(x)&=&\varphi*\rho^{(1)}_\eps(x),
\end{eqnarray*}
where $\rho^{(1)}_\eps$ [resp., $\rho^{(2)}_\eps$ and $\rho
^{(3)}_\eps$]
are the mollifiers in $\mR^{d+k+1}$ (resp., $\mR^{d+1}$ and~$\mR^d$). It
is clear that
\[
\|\nabla G^\eps_t\|_\infty+\|\nabla f^\eps_t\|_\infty+\|\nabla
\varphi
^\eps\|_\infty\leq K.
\]
By Theorem~\ref{Th3} and Proposition~\ref{Pro3},
there exists a time $T=T(K)$ such that for all $T\leq t\leq s\leq0$
and $x\in\mR^d$,
\[
\lim_{\eps\downarrow0}\mE|X^\eps_{t,s}(x)-X_{t,s}(x)|=0,
\]
where $X^\eps_{t,s}$ (resp., $X_{t,s}$) is the solution family of SFDE
(\ref{SFEE}) corresponding to
the coefficients ($G^\eps, f^\eps,\varphi^\eps$) [resp.,
($G,f,\varphi$)].
Using this limit, and by the dominated convergence theorem, it is easy
to verify that
for each $(t,x)\in[T,0]\times\mR^d$,
%
\begin{equation}\label{LL1}
u^\eps_t(x)\to u_t(x),
\end{equation}
where $u^\eps_t(x)$ is defined through $\varphi^\eps,f^\eps$ and
$X^\eps
_{t,s}(x)$ as in (\ref{EU}).
Moreover, by Proposition~\ref{Pro3}, we also have
%
\begin{equation}\label{EU1}
\sup_{\eps\in(0,1)}\sup_{t\in[T,0]}\|\nabla u^\eps_t\|_\infty
+\sup_{t\in[T,0]}\|\nabla u_t\|_\infty\leq C_{T,K}<+\infty.
\end{equation}
On the other hand, thanks to (\ref{EU1}), by Theorem~\ref{Th8}, there
exists another
time $T'=T'(K)\in[T,0)$ independent of $\eps$ such that
\[
u^\eps_t(x)=\varphi^\eps(x)+\int^0_t\cL_0 u^\eps_r(x)\,\dif r
+\int^0_tG^i_r(x,u^\eps_r(x))\partial _iu^\eps_r(x)\,\dif r+\int
^0_tf^\eps
_r(x)\,\dif r.
\]
In particular, for all $\psi\in C^\infty_0(\mR^d)$ and all $t\in[T',0]$,
%
\begin{eqnarray}
\langle u^\eps_t,\psi\rangle&=&\langle\varphi^\eps,\psi\rangle
+\int
^0_t\langle u^\eps_r,\cL_0^*\psi\rangle\,\dif r
+\int^0_t\langle G^i_r(u^\eps_r)\partial _iu^\eps_r,\psi\rangle
\,\dif
r
\nonumber
\\[-8pt]
\\[-8pt]
\nonumber
&&{}+\int
^0_t\langle f^\eps_r,\psi\rangle\,\dif r.
\end{eqnarray}
We want to take limits for both sides of the above identity by (\ref
{LL1}). The key point is to prove
\[
\int^0_t\langle G^{\eps,i}_r(u^\eps_r)\partial _iu^\eps_r,\psi
\rangle
\,\dif
r\to\int^0_t\langle G^i_r(u_r)\partial _iu_r,\psi\rangle\,\dif r,
\]
which will be obtained by proving the following two limits:
\begin{eqnarray*}
\int^0_t\bigl\langle\bigl(G^{\eps,i}_r(u^\eps_r)-G^i_r(u_r)\bigr)\partial
_iu^\eps
_r,\psi
\bigr\rangle\,\dif r&\to&0,\qquad \eps\to0,\\
\int^0_t\langle G^i_r(u_r)\partial _i(u^\eps_r-u_r),\psi\rangle
\,\dif r&\to&
0,\qquad \eps\to0.
\end{eqnarray*}
The first limit is clear by (\ref{LL1}), (\ref{EU1}) and the dominated
convergence theorem.
The second limit follows by (\ref{LL1}), (\ref{EU1}) and the
integration by parts formula.
\end{pf}

Now we are in a position to give:

\subsection{\texorpdfstring{Proof of Theorem \protect\ref{Main}}{Proof of Theorem 4.2}}\label{sec4.3}

We divide the proof into three steps.

(\textit{Step} 1). For $h\in\cB(\mR_-;\mW^{1,\infty}(\mR^d;\mR^k))$, define
\[
f^h_r(x):=F_r(x,h_r(x))
\]
and
\[
\sK:=\sup_{s\in\mR_-}(\|\nabla G_s\|_\infty+\|\nabla F_s\|
_\infty
)+\|\nabla\varphi\|_\infty.
\]
In this step, we prove the following claim:

\begin{claim*}\label{cl} For given $U\geq4\|\nabla\varphi\|_\infty$, there
exists a
time $T=T(\sK,U)<0$ such that
for any bounded measurable function $h\dvtx \mR_-\times\mR^d\to\mR^m$ satisfying
$\sup_{t\in[T,0]}\|\nabla h_t\|_\infty\leq U$,
it holds that
%
\begin{equation}\label{EPP2}
\|u^h_t\|_\infty\leq\|\varphi\|_\infty+|t|\sup_{s\in[t,0]}\|
f^h_s\|
_\infty,
\end{equation}
and
%
\begin{equation}\label{EPP1}
\sup_{t\in[T,0]}\|\nabla u^h_t\|_\infty\leq U,
\end{equation}
where $u^h_t(x)$ is defined by (\ref{EU}) in terms of $\varphi, f^h$
and $X^h_{t,s}(x)$,
and $\{X^h_{t,s}(x), T\leq t\leq s\leq0, x\in\mR^d\}$ is the unique
solution family of SFDE (\ref{SFEE})
corresponding to ($G,f^h,\varphi$).
\end{claim*}

\begin{pf*}{Proof of the \hyperref[cl]{Claim}} By Proposition~\ref{Pro3}, there exists a
time $T_1:=T_1(\sK, U)<0$
such that for all $x,y\in\mR^d$,
\[
\sup_{T_1\leq t\leq s\leq0}\mE|X_{t,s}(x)-X_{t,s}(y)|\leq2|x-y|.
\]
Using this and by the definition of $u^h_t(x)$ [see (\ref{EU})], we have
\[
|u^h_t(x)-u^h_t(y)|\leq2\|\nabla\varphi\|_\infty|x-y|+
2\int^0_t(\|\nabla_x F_r\|_\infty+\|\nabla_u F_r\|_\infty
U)|x-y|\,\dif r.
\]
So,
\begin{eqnarray*}
\sup_{s\in[t,0]}\|\nabla u^h_s\|_\infty&\leq&2\|\nabla\varphi\|
_\infty+
2|t|\sup_{s\in[t,0]}(\|\nabla_x F_s\|_\infty+\|\nabla_u F_s\|
_\infty
U)\\
&\leq&2\|\nabla\varphi\|_\infty+2|t|\sK(U+1).
\end{eqnarray*}
Since $U\geq4\|\nabla\varphi\|_\infty$,
choosing $T=\frac{U-2\|\nabla\varphi\|_\infty}{2\sK(U+1)}\wedge T_1$,
we obtain (\ref{EPP1}). Estimate (\ref{EPP2})
follows from definition (\ref{EU}).
\end{pf*}

(\textit{Step} 2). Set $u^0_t(x):=\varphi(x)$. We construct the following
iteration approximation sequence:
for $n\in\mN$,
\begin{eqnarray*}
X^n_{t,s}(x)&:=&X^{u^{n-1}}_{t,s}(x), \qquad u^n_t(x):=u^{u^{n-1}}_t(x),\\
f^n_t(x)&:=&f^{u^{n-1}}_t(x):=F_t(x,u^{n-1}_t(x)).
\end{eqnarray*}
By the above claim, there exists a time $T_1=T_1(\sK)<0$ such that for
all $n\in\mN$,
%
\begin{equation}\label{EPo1}
\qquad \|u^n_t\|_\infty\leq\|\varphi\|_\infty+|t|\sup_{s\in[t,0]}\|F_s\|
_\infty
, \qquad \sup_{t\in[T_1,0]}\|\nabla u^n_t\|_\infty\leq4\|\nabla\varphi\|
_\infty
.
\end{equation}
Hence,
\[
\|\nabla f^n_t\|_\infty\leq\|\nabla_x F_t\|_\infty+\|\nabla_u F_t\|
_\infty\|\nabla u^{n-1}_t\|_\infty
\leq\|\nabla_x F_t\|_\infty+4\|\nabla_u F_t\|_\infty\|\nabla
\varphi\|
_\infty.
\]
Thus, by the definition of $u^n_t(x)$ [see (\ref{EU})] and Proposition
\ref{Pro3} again,
there exists another time $T=T(\sK)\in[T_1,0)$ such that for all
$n,m\in
\mN$ and $t\in[T,0)$,
\begin{eqnarray*}
\|u^n_t-u^m_t\|_\infty&\leq&\|\nabla\varphi\|_\infty\sup_{x\in\mR
^d}\mE
|X^n_{t,0}(x)-X^m_{t,0}(x)|
+\int^0_t\|f^n_r-f^m_r\|_\infty\,\dif r\\
&&{} +\int^0_t\|\nabla f^n_r\|_\infty\sup_{x\in\mR^d}\mE
|X^n_{t,r}(x)-X^m_{t,r}(x)|\,\dif r\\
&\leq& C\int^0_t\|f^n_r-f^m_r\|_\infty\,\dif r\leq C\int^0_t\|
u^{n-1}_r-u^{m-1}_r\|_\infty\,\dif r,
\end{eqnarray*}
where $C$ is independent of $n,m$. By Gronwall's inequality, we obtain that
\[
\lim_{n,m\to\infty}\sup_{t\in[T,0]}\|u^n_t-u^m_t\|_\infty=0.
\]
Hence, there exists a $u_t\in\cB([T,0]\times\mR^d;\mR^k)$ such that
%
\begin{equation}\label{ER1}
\lim_{n\to\infty}\sup_{t\in[T,0]}\|u^n_t-u_t\|_\infty=0,
\end{equation}
and by (\ref{EPo1}),
\[
\|u_t\|_\infty\leq\|\varphi\|_\infty+|t|\sup_{s\in[t,0]}\|F_s\|
_\infty
,\qquad  \sup_{t\in[T,0]}\|\nabla u_t\|_\infty\leq4\|\nabla\varphi\|
_\infty.
\]
On the other hand, by Theorem~\ref{Th7}, $u^n_t(x)$ satisfies that for
all $\psi\in C^\infty_0(\mR^d;\mR^k)$,
\[
\langle u^n_t,\psi\rangle=\langle\varphi,\psi\rangle+\int
^0_t\langle
u^n_s,\cL_0^*\psi\rangle\,\dif r
+\int^0_t\langle G^i_s(u^n_s)\partial _iu^n_s+F_s(u^{n-1}_s),\psi
\rangle
\,\dif s.
\]
Thus, one can take limits as in Theorem~\ref{Th7} to obtain the
existence of a short time weak solution
for equation (\ref{EE1}). Moreover, (ii) follows from (\ref{EU}). The
existence of a maximal weak solution can be obtained as in the proof of
Theorem~\ref{Th2}
by shifting the time and the induction.
Thus, we conclude the proof of (i). As for (ii), it follows by (\ref{ER1})
and the definition of $u^n_t(x)$.

(\textit{Step} 3). Let $u\in\cB_{\mathrm{loc}}((T,0];\mW^{1,\infty}(\mR^d;\mR^k))$
be a
maximal weak solution of
equation (\ref{EE1}). Define for $(t,x)\in(T,0]\times\mR^d$,
\[
b_t(x):=G_t(x,u_t(x)), \qquad f_t(x):=F_t(x,u_t(x)).
\]
Then it is clear that
\[
b\in\cB_{\mathrm{loc}}((T,0];\mW^{1,\infty}(\mR^d;\mR^d)),\qquad  f\in\cB_{\mathrm{loc}}((T,0];\mW^{1,\infty}(\mR^d;\mR^k)).
\]
For $t\in(T,0]$, let $\{X_{t,s}(x),t\leq s\leq0,x\in\mR^d\}$ solve the
following SDE:
\[
X_{t,s}(x)=x+\int^s_tb_r(X_{t,r}(x))\,\dif r+\int^s_t\,\dif L_r, \qquad s\in[t,0].
\]
Define
%
\begin{equation}\label{EI1}
\tilde u_t(x):=\mE(\varphi(X_{t,0}(x)))+\int^0_t\mE
(f_s(X_{t,s}(x)))\,\dif s.
\end{equation}
By Theorem~\ref{Th6}, we have
\[
\tilde u_t(x)=u_t(x) \qquad \forall(t,x)\in(T,0]\times\mR^d.
\]
Suppose now that $T>-\infty$. For completing the proof, by (\ref{Max})
it is enough to show that
\[
\lim_{t\downarrow T}\|\nabla\tilde u_t(x)\|_\infty<+\infty.
\]
It immediately follows from (\ref{EI1}) and the following claim proved
in~\cite{Zh2}, Theorem~4.5,
which is stated in a slight variant.

\begin{claim*}
Under (\ref{EU3}) or $A$ nondegenerate, for any bounded
continuous function $\varphi$ and $T<t<s\leq0$,
\[
\|\nabla\mE\varphi(X_{t,s}(\cdot))\|_\infty\leq C_1(|t-s|\wedge
1)^{-1/\alpha}\|\varphi\|_\infty,
\]
where $C_1$ only depends on $d,\alpha,T$ and the bound of $b$.
\end{claim*}

\section*{Acknowledgment}
Deep thanks go to the referee for his/her very careful reading the
manuscript and very useful suggestions.


%


\printaddresses

\end{document}